\documentclass[preprint]{elsarticle}

\usepackage{lineno,hyperref}
\modulolinenumbers[5]

\journal{Journal of \LaTeX\ Templates}
\usepackage[cmex10]{amsmath}
\usepackage{amssymb}
\usepackage{algorithm}
\usepackage{algorithmic}

\newtheorem{theorem}{Theorem}

\usepackage{epstopdf}

\usepackage{subfig}
\usepackage{color}

\usepackage{tikz}
\usetikzlibrary{shapes,arrows,fit,positioning,shadows,calc}
\usetikzlibrary{plotmarks}
\usetikzlibrary{decorations.pathreplacing}
\usetikzlibrary{patterns}
\usetikzlibrary{automata}
\usepackage{pgfplots}
\pgfplotsset{compat=newest}
\usepackage{adjustbox}
\usepackage{caption}
\usepackage[a4paper,
            bindingoffset=0.2in,
            left =.5in,
            right=.5in,
            top  =1in,
            bottom=1in,
            footskip=.25in]{geometry}

\begin{document}

\begin{frontmatter}

\title{An Efficient Randomized QLP Algorithm for Approximating the Singular Value Decomposition}

\author[1]{M. F. Kaloorazi}
\ead{kaloorazi@xsyu.edu.cn}
\author[2]{K. Liu}
\ead{liukaide@mail.nwpu.edu.cn}
\author[2]{J. Chen}
\ead{jie.chen@nwpu.edu.cn}
\author[3]{R. C. de Lamare}
\ead{delamare@cetuc.puc-rio.br}


\address[1]{School of Electronic Engineering, Xi’an Shiyou University, Xi'an, China}
%addressline={Radarweg 29},
%postcode={1043 NX},
%city={Xi'an},
%country={China}}

\address[2]{CIAIC, Northwestern Polytechnical University, Xi'an, China}
%addressline={Radarweg 29},
%postcode={1043 NX},
%city={Xi'an},
%country={China}}

\address[3]{CETUC, Pontifical Catholic University of Rio de Janeiro, Rio de Janeiro, Brazil}
%addressline={Radarweg 29},
%postcode={1043 NX},
%city={Rio de Janeiro},
%country={Brazil}}

\begin{abstract}
In this paper, we introduce a randomized QLP decomposition called Rand-QLP. Operating on a matrix $\bf A$, Rand-QLP gives ${\bf A}={\bf QLP}^T$, where $\bf Q$ and $\bf P$ are orthonormal, and $\bf L$ is lower-triangular. Under the assumption that the rank of the input matrix is $k$, we derive several error bounds for Rand-QLP: bounds for the first $k$ approximate singular values and for the trailing block of the middle factor $\bf L$, which show that the decomposition is rank-revealing; bounds for the distance between approximate subspaces and the exact ones for all four fundamental subspaces of a given matrix; and bounds for the errors of low-rank approximations constructed by the columns of $\bf Q$ and $\bf P$. 
Rand-QLP is able to effectively leverage modern computational architectures, due to the utilization of random sampling and the \emph{unpivoted} QR decomposition, thus addressing a serious bottleneck associated with classical algorithms such as the singular value decomposition (SVD), column-pivoted QR (CPQR) and most recent matrix decomposition algorithms. To assess the performance behavior of different algorithms, we use an Intel Xeon Gold 6240 CPU running at 2.6 GHz with a NVIDIA GeForce RTX 2080Ti GPU. In comparison to CPQR and the SVD, Rand-QLP respectively achieves a speedup of up to 5 times and 6.6 times on the CPU and up to 3.8 times and 4.4 times with the hybrid GPU architecture. In terms of quality of approximation, our results on synthetic and real data show that the approximations by Rand-QLP are comparable to those of pivoted QLP and the optimal SVD, and in most cases are considerably better than those of CPQR.     
\end{abstract}

\begin{keyword}
Singular value decomposition, pivoted QLP, rank-revealing decomposition, randomized sampling, high-performance computing, hybrid architecture, communication-avoiding algorithm.
\end{keyword}

\end{frontmatter}

%\linenumbers
%{\footnotesize \textsuperscript{*} The first version of this paper appeared in arXiv in 2021 (arXiv:2110.01011)}

\section{Introduction}
Matrix decomposition is a factorization of a matrix into simpler constituents. Matrix decompositions are powerful computational tools for solving a wide array of problems that appear in signal processing and machine learning applications. Their use is due to the fact that these applications require information about orthonormal bases for the row and/or column spaces of a matrix (the bases that span the four fundamental subspaces of the matrix), or the eigen/singular-vales, the numerical rank, the condition number, or a low-rank representation of a matrix or a system. Applications of matrix decompositions include principal component analysis \cite{FanSZZ18}, linear least squares problems \cite{Bjorck15}, recommender systems \cite{LedentAK21}, image deconvolution and reconstruction \cite{GarveyMN18, MFKJCTSP20}, molecular dynamics simulations \cite{Gygi08}, coarse quantization in wireless networks \cite{DanaeedeN22}, community detection \cite{MaSZ21}, model matching for model order reduction \cite{Freund03}, background modeling \cite{MFKDeTSP18}, vibrations of fluid-solid structures \cite{Unger13}, subspace estimation over networks \cite{ChenRS17}, video compression \cite{CheW22}, seismic facies analysis for oil and gas reservoirs characterization \cite{RoyDM13}, meta-learning linear regression \cite{pmlr-v119-kong20a}, and Gaussian processes \cite{LiuHOSC20}.

The singular value decomposition (SVD) \cite{GolubVanLoan96}, the column-pivoted QR (CPQR), or rank-revealing QR factorization \cite{Chan87}, the  rank-revealing UTV decompositions \cite{StewartURV92, StewartULV93, Hansen98}, and the pivoted QLP decomposition \cite{StewartQLP} are widely used to factorize a general matrix. A major bottleneck associated with the computation of these factorizations is the communication cost \cite{ButtariLKD08, Hoemmen10, DemGGX15}. This cost is prohibitive, makes execution of these algorithms on modern computers challenging and in turn stymies their application to large matrices. The communication cost involves moving data between different levels of the memory hierarchy or/and between processors working in parallel. It is associated with the exploitation of level-1, 2 and 3 BLAS (the Basic Linear Algebra Subprograms) routines \cite{LawsonHKK79, DongarraDHH88, DongarraDHD90, Dongarraetal18, MFKJCTSPPbP21} by any algorithm: level-1 and level-2 BLAS routines are memory-bound and can not attain high performance on advanced architectures. Level-3 BLAS routines are, however, CPU-bound; they can harness the data locality in multicore architectures, thus attaining higher performance, very close to that rendered by processors. Operations executed in level-3 BLAS can be several times faster than those executed in level-1 and level-2 BLAS.

%\noindent \textbf{Our Contributions.}
\subsection{Contributions}
This paper develops a simple, numerically stable, yet fast and readily parallelized decomposition algorithm that, like the pivoted QLP and the SVD, provides all information on singular values and four fundamental subspaces of a given matrix.  The algorithm makes use of randomization and the decomposition takes the form of QLP and is called Rand-QLP. 
The computation of Rand-QLP requires only matrix-matrix multiplication and the highly parallelizable unpivoted QR decomposition. Hence the operations of Rand-QLP can be cast almost entirely in level-3 BLAS. 

We present a theoretical error analysis for Rand-QLP. We show that Rand-QLP is rank-revealing by deriving bounds for the leading approximate singular values as well as the trailing block of the middle factor. We furnish bounds for the distance between approximate subspaces and the exact ones for all four fundamental subspaces of a given matrix. We further bound from above the errors of low-rank approximations constructed by Rand-QLP. All of our bounds take into account the distribution of the random matrix used.

We make use of a hybrid GPU-accelerated multicore system to empirically compare the performance of Rand-QLP and several classical and recently developed algorithms on randomly generated data as well as real data.

\subsection{Organization}
The organization of this paper is as follows. Section \ref{secRelWor} reviews related classical and recent works concerning matrix decomposition algorithms. Section \ref{secProAlg} describes the proposed Rand-QLP algorithm and discusses the intuition behind it. Section \ref{secErrBou} presents the theoretical error bounds for Rand-QLP, whose proofs are given in the appendix, and details the computational cost of the algorithm. Section \ref{secNuSim} presents and discusses simulation results of our numerical tests, and concluding remarks are given in Section \ref{secCon}.

\section{Related Work}
\label{secRelWor}
Let ${\bf A} \in \mathbb R^{m \times n}$, with $m \ge n$, and its rank be $k$. The SVD is considered the best decomposition of a matrix; it provides information on singular values and four fundamental subspaces of the matrix. The (thin) SVD \cite[Section 2.5.4]{GolubVanLoan96} of $\bf A$ is written as: 
\begin{equation}
{\bf A} = {\bf U}{\bf \Sigma}{\bf V}^T= \begin{bmatrix} {{\bf U}_k \quad {\bf U}_\perp} \end{bmatrix}
  \begin{bmatrix}
       {\bf \Sigma}_k & {\bf 0}  \\
       {\bf 0} & {\bf \Sigma}_\perp
  \end{bmatrix}\begin{bmatrix}{{\bf V}_k \quad {\bf V}_\perp} \end{bmatrix}^T,
\label{eqSVD}
\end{equation}
where ${\bf U}_k \in \mathbb R^{m \times k}$ and ${\bf U}_\perp \in \mathbb R^{m \times (n-k)}$ are orthonormal and span $\mathcal{R}({\bf A})$ and $\mathcal{N}({\bf A}^T)$, respectively. The diagonal ${\bf \Sigma} \in \mathbb R^{n \times n}$ contains the singular values $\sigma_i$s in a non-increasing order: $\sigma_1\ge\sigma_2\ge ... \ge \sigma_n$, where ${\bf \Sigma}_k$ comprises the first $k$ and ${\bf \Sigma}_\perp$ the remaining $n-k$ singular values. ${\bf V}_k \in \mathbb R^{n \times k}$ and ${\bf V}_\perp \in \mathbb R^{n \times (n-k)}$ are orthonormal and span $\mathcal{R}({\bf A}^T)$ and $\mathcal{N}({\bf A})$, respectively. ($\mathcal{R}(\cdot)$ and $\mathcal{N}({\cdot})$ indicate the range and null space.) The SVD can be used to determine the rank of a matrix, to  compute angles between subspaces as well as optimal low-rank approximations.

The unpivoted QR decomposition \cite[Section 5.2]{GolubVanLoan96} gives ${\bf A} = {\bf Q}_\text{A}{\bf R}_\text{A}$, where ${\bf Q}_\text{A} \in \mathbb R^{m \times n}$ is orthonormal and contains the bases for the columns of $\bf A$, and ${\bf R}_\text{A} \in \mathbb R^{n \times n}$ is upper triangular. (This is in fact the ``reduced" version of the decomposition, for silent columns and rows in the two factors were discarded \cite[Lecture 7]{TrefethenBau97}.) The QR decomposition is computed through the Gram–Schmidt process \cite[Section 2.3.4]{Bjorck15}, \cite[Lecture 8]{TrefethenBau97} or Householder reflections \cite[Section 5.1.2]{GolubVanLoan96}, \cite[Lecture 10]{TrefethenBau97}. A variant of the QR decomposition is the column-pivoted QR (CPQR) \cite[Section 5.4.1]{GolubVanLoan96}, or rank-revealing QR decomposition \cite{Chan87}, \cite[Section 2.4.3]{Bjorck15}. It is computed by Householder orthogonal triangularization where columns with the largest magnitudes are swapped with other columns before the reduction proceeds. CPQR factorizes $\bf A$ as:
\begin{equation}
{\bf A}{\bf \Pi} = {\bf Q}_\text{piv}{\bf R}_\text{piv}= {\bf Q}_\text{piv}
  \begin{bmatrix}
       {\bf R}_{11} & {\bf R}_{12}  \\
       {\bf 0} & {\bf R}_{22}
  \end{bmatrix},
\label{equCPQR}
\end{equation}
where ${\bf \Pi}\in \mathbb R^{n \times n}$ is a permutation matrix. CPQR is rank-revealing in the sense that the following two conditions hold for ${\bf R}_{11} \in \mathbb R^{k \times k}$ and ${\bf R}_{22} \in \mathbb R^{(n-k) \times (n-k)}$:
\begin{equation}\label{eqDeterRankR}
\begin{aligned}
\sigma_\text{min}({\bf R}_{11}) \ge \sigma_k({\bf A})/f(n, k),\quad \text{and} \quad
\|{\bf R}_{22}\|_2 \le \sigma_{k+1}({\bf A})g(n,k),
\end{aligned}
\end{equation}
where $f(n,k)$ and $g(n,k)$ are low degree polynomials in $n$ and $k$, and $\|\cdot\|_2$ indicates the $\ell_2$-norm. Elements on the diagonal of ${\bf R}_\text{piv}$ are approximations to the singular values of $\bf A$. CPQR, however, does not explicitly establish orthonormal bases for the rows of $\bf A$.

Introduced by Stewart in \cite{StewartURV92} and \cite{StewartULV93}, the rank-revealing URV and ULV decompositions are a compromise between the SVD and CPQR; they are more efficient than the SVD in terms of arithmetic cost, and they furnish orthonormal bases for $\mathcal{R}({\bf A}^T)$ and $\mathcal{N}({\bf A})$. These algorithms, often called UTV decompositions, factor $\bf A$ such as ${\bf A}= {\bf UTV}^T$; ${\bf U} \in \mathbb R^{m \times n}$ and ${\bf V} \in \mathbb R^{n \times n}$ are orthonormal, and ${\bf T} \in \mathbb R^{n \times n}$ is either upper triangular (the URV) or lower triangular (the ULV). The computation of UTV decompositions involves two main steps: an initial QR decomposition, and deflation steps to estimate the singular values \cite{FierroHan97, Hansen98}.

The pivoted QLP (p-QLP) decomposition \cite{StewartQLP} is another rank-revealing algorithm which is built on CPQR. Specifically, let $\bf A$ have a CPQR such that ${\bf A\Pi} = {\bf Q}_\text{piv}{\bf R}_\text{piv}$. Then a CPQR is performed on ${\bf R}_\text{piv}^T$ such that ${\bf R}_\text{piv}^T\acute{\bf \Pi} = \acute{\bf Q}\acute{\bf R}$. Hence p-QLP is given by:
\begin{equation}\notag
{\bf A} = {\bf Q}_\text{piv}\acute{\bf \Pi}\acute{\bf R}^T\acute{\bf Q}^T{\bf \Pi}^T  \triangleq 
{\bf Q}_\text{qlp}{\bf L}_\text{qlp}{\bf P}_\text{qlp}^T,
\end{equation}
where ${\bf Q}_\text{qlp}\triangleq {\bf Q}_\text{piv}\acute{\bf \Pi}$ and ${\bf P}_\text{qlp}\triangleq {\bf \Pi}\acute{\bf Q}$ are orthonormal and contain bases for the columns and rows of $\bf A$, respectively, and the diagonal entries of ${\bf L}_\text{qlp} \triangleq \acute{\bf R}^T$ are approximations to the singular values of $\bf A$. (For other matrix decompositions see \cite{GolubVanLoan96,Bjorck15}.)

In order to compute the SVD of a matrix, the classical bidiagonalization method is usually used \cite{GolubVanLoan96, Dongarraetal18}, which comprises two main steps: reduction of the matrix to a bidiagonal form, and reduction of the bidiagonal form to a diagonal form. Most operations of these two steps are in level-1 and level-2 BLAS. Due to employment of the pivoting strategy, roughly half of the operations of CPQR are in BLAS-3 \cite{QuintanaSB98}. The UTV algorithm in the second step deflates large or small singular values one at a time \cite{FierroHan97, Hansen98}. A large portion of these operations are in level-1 and level-2 BLAS. Most operations of the unpivoted QR decomposition are, however, in level-3 BLAS \cite{ButtariLKD08, Hoemmen10}.

Recently, Nakatsukasa and Higham \cite{NakatsukasaH13} developed the QR-based Dynamically Weighted Halley (QDWH)-SVD algorithm, a communication friendly algorithm for computing the SVD of a matrix. The building block of the algorithm is the computation of the polar decomposition \cite[Chapter 8]{Higham08}, which is done by a Halley iteration using the QR decomposition. QDWH-SVD consists of two main steps: i) computation of the polar decomposition ${\bf A}={\bf U}_p{\bf H}$, where ${\bf U}_p \in \mathbb R^{m\times n}$ is orthonormal, and ${\bf H} \in \mathbb R^{n\times n}$ is Hermitian positive semidefinite, and ii) computation of the eigendecomposition ${\bf H}={\bf V}{\bf \Sigma}{\bf V}^T$. Hence ${\bf A}=({\bf U}_p{\bf V}){\bf \Sigma}{\bf V}^T$. Martinsson, Quintana-Orti and Heavner \cite{MartinssonQH19} presented randUTV, a randomized algorithm that approximates the UTV decomposition of a matrix. Using randomization and the SVD, combined with power iteration technique and ``oversampling", randUTV factors one block of b columns of the matrix at a time. By piecing the results together, it constructs a UTV-like factorization of the matrix. A more detailed comparison of QDWH-SVD, randUTV and Rand-QLP in computational procedure, flop count and performance is given in Section \ref{subsRuntime}.

\section{Proposed Rand-QLP Algorithm}
\label{secProAlg}

This section describes our Rand-QLP algorithm. The mechanics of the algorithm are very simple and involve only matrix-matrix multiplication and the QR decomposition. These operations can be cast almost completely in terms of level-3 BLAS, which enables Rand-QLP to harness the power of modern computational environments. In addition, for the singular values and any rank-$r$ approximations, where $1\le r \le n$, the accuracy of the approximations of Rand-QLP are comparable to those of the optimal SVD. These properties make Rand-QLP preferable to classical algorithms such as p-QLP and the SVD. 
 
\subsection{Description}
The algorithmic procedure of Rand-QLP to factorize the matrix $\bf A$ defined above is as follows. First, a random matrix ${\bf \Omega} \in \mathbb R^{m \times n}$ (in this work standard Gaussian) is generated and left-multiplied by ${\bf A}^T$. Next, orthonormal bases for the columns of  ${\bf A}^T{\bf \Omega}$ are obtained, e.g., via the QR decomposition:  
\begin{equation}\notag
\bar{\bf Q} = \texttt{orth}({\bf A}^T{\bf \Omega}) \quad \text{where} \quad \bar{\bf Q} \in \mathbb R^{n \times n}.
\end{equation}
Then $\bar{\bf Q} \in \mathbb R^{n \times n}$ is left-multiplied by ${\bf A}$, and orthonormal bases for the columns of the product ${\bf A}\bar{\bf Q}$ are computed:
\begin{equation}\notag
{\bf Q} = \texttt{orth}({\bf A}\bar{\bf Q}) \quad \text{where} \quad {\bf Q} \in \mathbb R^{m \times n}.
\end{equation}
We then form ${\bf Q}^T{\bf A}$ and compute the QR decomposition of its transpose:
\begin{equation}\notag
{\bf PR} =({\bf Q}^T{\bf A})^T \quad \text{where} \quad {\bf P, R} \in \mathbb R^{n \times n}.
\end{equation}
Thus ${\bf A}={\bf QLP}^T$, where ${\bf L} \triangleq {\bf R}^T$. Algorithm \ref{Alg1} presents the
procedure for computing Rand-QLP. 

\begin{algorithm}
\caption{Rand-QLP}
\renewcommand{\algorithmicrequire}{\textbf{Input:}}
\begin{algorithmic}[1]
\REQUIRE ~ % ：Input
An $m \times n$ matrix ${\bf A}$.
\renewcommand{\algorithmicrequire}{\textbf{Output:}}
\REQUIRE  ${\bf Q} \in \mathbb R^{m \times n}$ is orthonormal and contains the approximate bases for $\mathcal{R}({\bf A})$ and $\mathcal{N}({\bf A}^T)$; ${\bf L}\in \mathbb R^{n \times n}$ is lower triangular and its diagonals approximate the singular values; ${\bf P} \in \mathbb R^{n \times n}$ is orthonormal and contains the approximate bases for $\mathcal{R}({\bf A}^T)$ and $\mathcal{N}({\bf A})$.
\STATE \textbf{function} {[${\bf Q}, {\bf L}, {\bf P}$]}{}$=$\texttt{{Rand\_QLP}}(${\bf A}$)
\STATE \quad ${\bf \Omega}=\texttt{randn}(m, n)$
\STATE \quad $\bar{\bf Q} = \texttt{orth}({\bf A}^T{\bf \Omega})$ 
\STATE \quad ${\bf Q} = \texttt{orth}({\bf A}\bar{\bf Q})$
\STATE \quad ${\bf PR} =({\bf Q}^T{\bf A})^T \rightarrow {\bf L} \triangleq {\bf R}^T$ 
\STATE \textbf{end function}
\end{algorithmic}\label{Alg1}
\end{algorithm}

\subsection{Intuition} 
To motivate the algorithm and to establish its error analysis, we first elaborate on the computation of p-QLP. Two CPQR decompositions are required to form p-QLP. However, as Stewart argues \cite{StewartQLP}, the second CPQR is not necessary, due to column pivoting done by the first CPQR. In other words, if in the first step of the p-QLP computation a matrix ${\bf Q}_\text{piv}$ is formed such that whose first $k$ columns span $\mathcal{R}({\bf A})$ and whose last $n-k$ columns span $\mathcal{N}({\bf A}^T)$, as in \eqref{equCPQR}, the second CPQR performed on ${\bf R}_\text{piv}^T$ can then be supplanted by an unpivoted QR decomposition. In the computation of Rand-QLP, through right multiplication of ${\bf A}^T$ by ${\bf \Omega}$ and of ${\bf A}$ by $\bar{\bf Q}$, and an orthonormalization thereafter, we construct a matrix $\bf Q$ that gives a good approximation to ${\bf Q}_\text{piv}$ in \eqref{equCPQR}. The matrix ${\bf Q}^T{\bf A}$ is in fact an approximation to ${\bf R}_\text{piv}$, and therefore we apply the QR decomposition to its transpose.    

%%One way to convert Rand-QLP to an SVD-like decomposition is through first block diagonalization of the middle factor $\bf L$. Under the assumption that there is a well-defined gap in the spectrum separating the first $k$ singular values from the remaining ones, there is a lower-triangular matrix ${\bf W}\triangleq \begin{bmatrix}
%{\bf I} & {\bf 0}  \\
%\bar{\bf W} & {\bf I}  
%\end{bmatrix}$, whose inverse is ${\bf W}^{-1}= \begin{bmatrix}
%{\bf I} & {\bf 0}  \\
%-\bar{\bf W} & {\bf I}  
%\end{bmatrix}$, 
%where $\bar{\bf W} \in \mathbb R^{n-k \times k}$, such that  
%\begin{equation}%\notag
%{\bf W}^{-1}{\bf LW}=
%\begin{bmatrix}
%{\bf L}_{11} & {\bf 0}  \\
%{\bf 0} & {\bf L}_{22}  
%\end{bmatrix}.
%\end{equation}
%This is obtained by solving the \emph{Sylvester equation} of the form: 
%\begin{equation}
%{\bf L}_{22}\bar{\bf W} - \bar{\bf W}{\bf L}_{11} = -{\bf L}_{21}.
%\end{equation}
%The SVD-like deccomposition may then be obtained by factorizing the smaller blocks on the diagonal and piecing the components together.

\section{Error Bounds and Computational Cost}
\label{secErrBou}
This section presents our theoretical analysis as well as the computational complexity of Rand-QLP. Theorem \ref{ThSinRQLP} asserts that Rand-QLP is rank-revealing through bounding the singular values of the diagonal blocks of $\bf L$. Theorem \ref{ThSins} furnishes upper bounds for the distance between approximate subspaces and the exact ones for all four fundamental subspaces of an input matrix, while Theorem \ref{ThLowRankErr} bounds the errors of low-rank approximations constructed by Rand-QLP. We have used a random matrix with a standard Gaussian distribution for computing Rand-QLP and all bounds presented take into account this distribution.

\subsection{Error Bounds} 
Let Rand-QLP compute a factorization for $\bf A$ as follows:
\begin{equation}%\notag
{\bf A} = {\bf QLP}^T=[{\bf Q}_1\quad {\bf Q}_2]
\begin{bmatrix}
{\bf L}_{11} & {\bf 0}  \\
{\bf L}_{21} & {\bf L}_{22}  
\end{bmatrix}
\begin{bmatrix} 
{\bf P}_1^T\\
{\bf P}_2^T
\end{bmatrix},
\end{equation}
where ${\bf Q}_1 \in \mathbb R^{m \times k}$,  ${\bf Q}_2 \in \mathbb R^{m \times (n-k)}$, ${\bf L}_{11} \in \mathbb R^{k \times k}$,  ${\bf L}_{22} \in \mathbb R^{(n-k) \times (n-k)}$, ${\bf P}_1 \in \mathbb R^{n \times k}$, and ${\bf P}_2 \in \mathbb R^{n \times (n-k)}$. Let further ${\bf \Omega}_1 \in \mathbb R^{m \times k}$ be formed by the first $k$ columns of the random matrix $\bf \Omega$ used to compute Rand-QLP, and $\widetilde{\bf \Omega} \triangleq {\bf U}^T {\bf \Omega}_1  = [\widetilde{\bf \Omega}_{1}^T \quad \widetilde{\bf \Omega}_{2}^T]^T$, where $\widetilde{\bf \Omega}_1 \in \mathbb R^{k\times k}$ and $ \widetilde{\bf \Omega}_2 \in \mathbb R^{(n-k)\times k}$. (We recall that $\bf U$ contains $\bf A$'s left singular vectors.) 

\subsubsection{Singular values} The theorem below bounds the $k$ singular values of ${\bf L}_{11}$ and the largest singular value of ${\bf L}_{22}$, thus showing that Rand-QLP is a rank-revealer.  
\begin{theorem}\label{ThSinRQLP}
Let ${\bf A}$ be an ${m\times n}$ matrix whose SVD is defined in \eqref{eqSVD}, $\bf L$ be constructed by Rand-QLP whose diagonals are $\widehat{\sigma}_i$, and $\psi_i= \frac{\sigma_{k+1}}{\sigma_i}$. Then, for $i = 1,..., k$, we have
\begin{equation}\label{eqSinV}
\sigma_i \ge \widehat{\sigma}_i\ge \frac{\sigma_i}{\sqrt{1 + \psi_i^4\|\widetilde{\bf \Omega}_2\widetilde{\bf \Omega}_1^{-1}\|_2^2}}.
\end{equation}

\begin{equation}\label{eqL22B}
\begin{aligned}
\|{\bf L}_{22}\|_2\le \big(1 + \psi_k^2\|\widetilde{\bf \Omega}_{2}\widetilde{\bf \Omega}_{1}^{-1} \|_2\big)\sigma_{k+1}.
\end{aligned}
\end{equation}
\end{theorem}

\textit{Proof.} Proof of this theorem is presented in Appendix \ref{AppdxPrTh1}.

The results in Theorem \ref{ThSinRQLP} parallel those of the deterministic algorithms given in \eqref{eqDeterRankR}. However, in contrast, they take into account the randomness involved.

It is possible to bound all singular values computed by Rand-QLP, but rather with a different approach. This requires a result from \cite{MathiasStewart93} (see also \cite{HuckabyChan03}). Specifically, let $\widehat{\bf L}$ and $\tau$ be defined as follows:
\begin{equation}\notag
\widehat{\bf L}=\begin{bmatrix}
{\bf L}_{11} & {\bf 0}  \\
{\bf 0} & {\bf L}_{22}  
\end{bmatrix},
\quad \text{and} \quad \tau= \dfrac{\|{\bf L}_{22}\|_2}{\sigma_k({\bf L}_{11})}<1.
\end{equation}
Then for $j=1,..., n$,
\begin{equation}\notag
\dfrac{\sigma_j({\bf L})}{\sigma_j(\widehat{\bf L})}= 1+ O\Bigg(\dfrac{\|{\bf L}_{12}\|_2^2}{(1-\tau)^2\sigma_k({\bf L}_{11})^2}\Bigg).
\end{equation}

Huckaby and Chan \cite{HuckabyChan03} have used the above result for bounding all singular values of a matrix decomposed by p-QLP; see \cite[Theorem 3.4]{HuckabyChan03}. 

\subsubsection{Distance between subspaces}
The distance between subspaces \cite{GolubVanLoan96}, denoted here by dist$(\cdot, \cdot)$, is used to measure the closeness of two subspaces; it corresponds to the largest canonical angle between the subspaces. Let 
\begin{equation}
\begin{aligned}
& \text{sin}\theta_Q = \text{dist}(\mathcal{R}({\bf U}_k), \mathcal{R}({\bf Q}_1)). \\
& \text{sin}\phi_Q = \text{dist}(\mathcal{R}({\bf U}_\perp), \mathcal{R}({\bf Q}_2)).\\
& \text{sin}\theta_P = \text{dist}(\mathcal{R}({\bf V}_k), \mathcal{R}({\bf P}_1)).\\
&\text{sin}\phi_P = \text{dist}(\mathcal{R}({\bf V}_\perp), \mathcal{R}({\bf P}_2)).
\end{aligned} \notag
\end{equation}

The following theorem bounds the sine of canonical angles between all four estimated subspaces and the exact subspaces given by the SVD.  
\begin{theorem} \label{ThSins}
With the notation and hypotheses of Theorem \ref{ThSinRQLP}, if $\|{\bf L}_{22}\|_2< \sigma_k({\bf L}_{11})$ then
\begin{equation}
\begin{aligned}
& \text{sin}\theta_Q \le \psi_k^2\|\widetilde{\bf \Omega}_{2}\widetilde{\bf \Omega}_{1}^{-1}\|_2,  \hspace{2cm}
 \text{sin}\theta_P \le \psi_k^3\|\widetilde{\bf \Omega}_{2}\widetilde{\bf \Omega}_{1}^{-1}\|_2, \\ 
& \text{sin}\phi_Q \le  \frac{\psi_k^2\|\widetilde{\bf \Omega}_{21}\widetilde{\bf \Omega}_{11}^{-1} \|_2^2}{\sqrt{1+\psi_k^4\|\widetilde{\bf \Omega}_{2}\widetilde{\bf \Omega}_{1}^{-1} \|_2^2}}, \hspace{0.9cm}
 \text{sin}\phi_Q \le  \frac{\psi_k^3\|\widetilde{\bf \Omega}_{21}\widetilde{\bf \Omega}_{11}^{-1} \|_2^2}{\sqrt{1+\psi_k^6\|\widetilde{\bf \Omega}_{2}\widetilde{\bf \Omega}_{1}^{-1} \|_2^2}}.
\end{aligned}
\notag
\end{equation} 
\end{theorem}

\textit{Proof.} Proof of this theorem is given in Appendix \ref{AppdxPrTh2}.

\subsubsection{Rank-k approximation error}
The following theorem furnishes in terms of $\ell_2$- and Frobenius norm upper bounds for the errors of rank-$k$ approximations constructed by Rand-QLP. 

\begin{theorem} \label{ThLowRankErr}
With the notation and hypotheses of Theorem \ref{ThSinRQLP},  let further ${\bf Q}_k$ and ${\bf P}_k$ be formed by the first $k$ columns of ${\bf Q}$ and ${\bf P}$, respectively, and $\rho=\{2, F\}$. Then
\begin{equation}
\begin{aligned}
& \|({\bf I} - {\bf Q}_k{\bf Q}_k^T){\bf A}\|_\rho \le \|{\bf \Sigma}_\perp\|_\rho+ \|{\bf A}\|_\rho \psi_k^2\|\widetilde{\bf \Omega}_{2}\widetilde{\bf \Omega}_{1}^{-1}\|_2.\\
& \|{\bf A}({\bf I} - {\bf P}_k{\bf P}_k^T)\|_\rho \le \|{\bf \Sigma}_\perp\|_\rho + \|{\bf A}\|_\rho \psi_k^3\|\widetilde{\bf \Omega}_{2}\widetilde{\bf \Omega}_{1}^{-1}\|_2.
\end{aligned}
\end{equation}
\end{theorem}

\textit{Proof.} Proof of this theorem is given in Appendix \ref{AppdxPrTh3}.

\subsection{Computational Cost}    
The computation of Rand-QLP requires the following number of flops (floating-point operations): generating the random matrix $\bf \Omega$ needs $mn$ flops; forming ${\bf A}^T{\bf \Omega}$ needs $2mn^2-n^2$ flops, and orthonormalization of this matrix to obtain $\bar{\bf Q}$ via the QR factorization needs $5n^3/3$ flops; computing ${\bf A}\bar{\bf Q}$ needs $2mn^2-mn$ flops; forming $\bf Q$ needs $2mn^2-n^3/3$ flops; computing ${\bf Q}^T{\bf A}$ requires $2mn^2-n^2$ flops, and the QR factorization on this matrix requires $5n^3/3$ flops. Hence
\begin{equation}\notag
C_\text{Rand-QLP} = 8mn^2+3n^3-2n^2.
\end{equation}

The above cost involves computation of all approximate subspaces as well as singular vales of $\bf A$ by Rand-QLP. 

The cost of computing the SVD is comparable to that of Rand-QLP \cite[p. 254]{GolubVanLoan96}. Computing CPQR, however, requires $2mn^2+n^3+4(m^2n-mn^2)$ flops, which is around three times less that that of Rand-QLP. However, as pointed out earlier, the computation of the SVD requires considerable amount of level-1 and level-2 BLAS operations \cite{Dongarraetal18}, and only half of the flops of CPQR are in level-3 BLAS \cite{QuintanaSB98}. This makes implementation of these two algorithms challenging on advanced computers, as they can not be efficiently parallelized. However, Rand-QLP can be computed almost entirely using level-3 BLAS routines, hence harnessing the structure of modern architectures. This claim is substantiated in Section \ref{secNuSim} through implementation of different algorithms on a hybrid GPU-accelerated machine. 

Rand-QLP readily lends itself to more complicated implementations due to the development of the blocked and tiled methods for computing the QR decomposition of dense and block low-rank matrices, which are amenable to high degrees of parallelism \cite{ButtariLKD08, AprY22}.

\section{Numerical Simulations}
\label{secNuSim}
We conduct numerical tests to investigate the performance of Rand-QLP. The results are compared with those of several existing algorithms. Tests are done in Python and the architecture used for this study is an Intel Xeon Gold 6240 CPU running at 2.6 GHz, with 18 cores and 251 GB of
RAM, equipped with a NVIDIA GeForce RTX 2080Ti GPU. 

\subsection{Runtime Comparison}
\label{subsRuntime}
We compare the speed of Rand-QLP against those of four algorithms, namely the SVD, CPQR, QDWH-SVD \cite{NakatsukasaH13}, and randUTV \cite{MartinssonQH19}. (In our comparison, we have discarded p-QLP as its cost can be over twice of that of CPQR, due to extra CPQR and matrix-matrix multiplications; see Section \ref{secRelWor} for details.) Before presenting the runtime results, we discuss the last two algorithms in more detail.
\begin{itemize}
\item QDWH-SVD \cite{NakatsukasaH13}. To furnish the SVD for $\bf A$, this algorithm computes two polar decompositions, one for $\bf A$ and one for $\bf H$, by using a Halley iteration. It performs the QR decomposition on matrices of size $(m+n)\times n$ for $\bf A$ and of size $2n\times n$ for $\bf H$. The number of iterations can reach up to six, depending on the condition number \cite[p. 94]{TrefethenBau97} of $\bf A$. The flop count for QDWH-SVD ranges from $8mn^2+(27+2/3)n^3$ to $24mn^2+(28+1/3)n^3$.
\item randUTV \cite{MartinssonQH19}. This algorithm uses randomized sampling and the SVD and gives a UTV-like factorization incrementally. The authors utilize the power iteration scheme and oversampling to improve the approximation accuracy. They argue that the cost of randUTV is three times higher than that of CPQR if the power iteration factor $q$ is set to zero; if $q=1$ or $q=2$, the cost would be four times or five times higher than that of CPQR, respectively. If a matrix $\bf A$ is stored externally, randUTV needs two passes over $\bf A$ if $q=0$, and $2q+2$ passes over $\bf A$ otherwise.    
\end{itemize}

The flop counts for Rand-QLP and randUTV with $q=0$ are comparable. Further, the latter requires two passes over an externally stored $\bf A$, while the former needs three passes. However, randUTV with $q \ge 1$ requires more arithmetic operations than Rand-QLP, as it needs more passes over $\bf A$. Rand-QLP is more efficient than QDWH-SVD in computational cost. 

We apply the five algorithms for factorization of the following matrices.
\begin{itemize}
\item {Synthetic Data.} These dense, $n\times n$ matrices have random entries from a uniform distribution on $(0,1)$. We note that for this task, the distribution of singular values is immaterial. 
\item {Real Data.} These data involve five sparse, square, non-symmetric matrices from the SuiteSparse Matrix Collection (formerly called the University of Florida Sparse Matrix Collection) \cite{DavisH11}. Information of these matrices are presented in Table \ref{TableMatSuiSpa}.
\end{itemize}

\begin{table}[!htb] 
\small{
\begin{center}   
\caption{Matrices from SuiteSparse Matrix Collection.}
\vspace{-.5cm}
\begin{tabular}
{p{1.2cm} p{1.5cm} p{6cm}}
\noindent\rule{9.2cm}{0.4pt}\\
Matrix & Dimensions & General problem category  \\
\noindent\rule{9.2cm}{0.4pt}\\
\texttt{watt1}   & $\quad$ 1856$^2$ & Computational fluid dynamics problem. \\
\texttt{meg1}    & $\quad$ 2904$^2$ & Circuit simulation problem. \\
\texttt{gemat11} & $\quad$ 4929$^2$ & Power network problem sequence. \\
\texttt{bayer03} & $\quad$ 6747$^2$ & Chemical process simulation Problem.  \\
\texttt{goodwin} & $\quad$ 7320$^2$ & Computational fluid dynamics problem. \\
\noindent\rule{9.2cm}{0.4pt}
\end{tabular}
\label{TableMatSuiSpa}
\end{center}}   
\end{table}

\begin{figure}[t]
\begin{center}       
\input{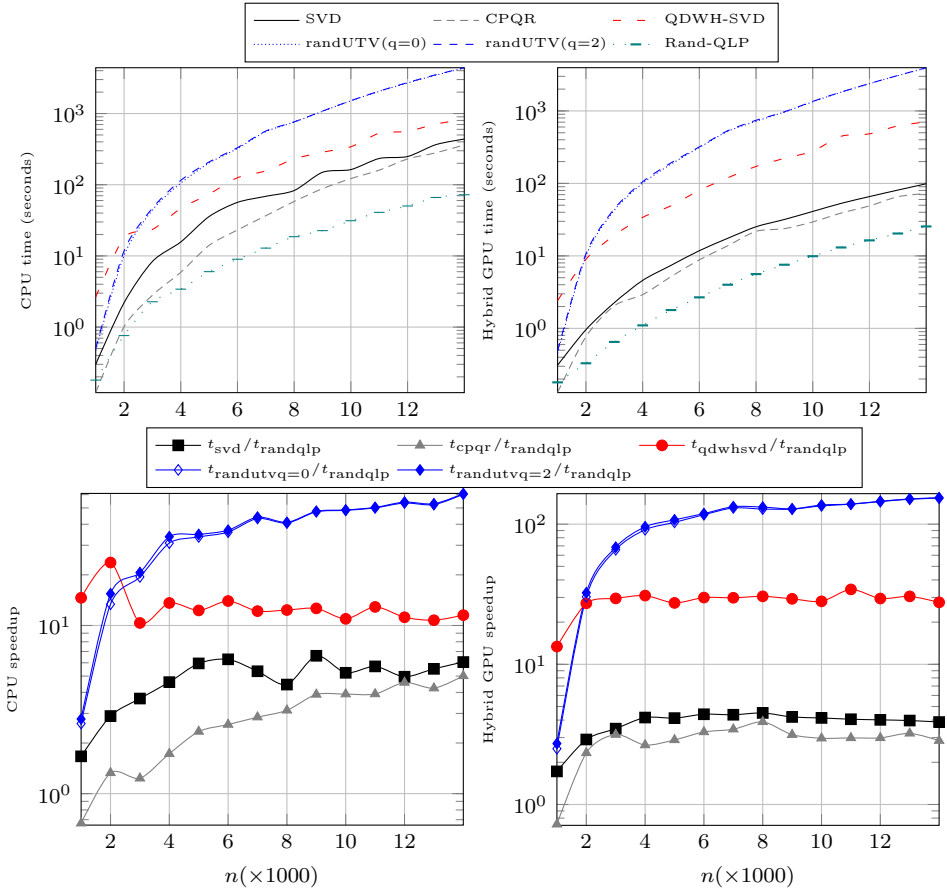}     
\captionsetup{justification=centering,font=scriptsize}  
\caption{Computational times for different algorithms on random matrices (top row) and speedups offered by Rand-QLP (bottom row). Our implementation of randUTV is that presented in \cite[Figure 3]{MartinssonQH19}, and we have set the block size $b$=100 and oversampling parameter to 5. (Note that the y-axis of the speedups are in logarithmic scale; the speedups of Rand-QLP over QDWH-SVD, and randUTV were so substantial that a meaningful plot demanded doing so.)} 
\label{figSpRandom}      
\end{center}
\end{figure}

\begin{figure}[t]
\begin{center}       
	% This file was created by matlab2tikz v0.4.7 running on MATLAB 8.3.
% Copyright (c) 2008--2014, Nico Schlömer <nico.schloemer@gmail.com>
% All rights reserved.
% Minimal pgfplots version: 1.3
% 
% The latest updates can be retrieved from
%   http://www.mathworks.com/matlabcentral/fileexchange/22022-matlab2tikz
% where you can also make suggestions and rate matlab2tikz.
% 
%
% defining custom colors
\usetikzlibrary{positioning,calc}

\definecolor{mycolor1}{rgb}{0.00000,1.00000,1.00000}%
\definecolor{mycolor2}{rgb}{1.00000,0.00000,1.00000}%

\pgfplotsset{every axis label/.append style={font=\footnotesize},
every tick label/.append style={font=\footnotesize}
}

\begin{tikzpicture}[font=\footnotesize] 

\begin{axis}[%
name=ber,
ymode=log,
width  = 0.32\columnwidth, 
height = 0.24\columnwidth, 
scale only axis,
xmin  = 1,
xmax  = 5,
xmajorgrids,
ymin = 0.57,
ymax = 593,
xtick       ={1, 2, 3, 4, 5},
xticklabels ={\tiny{watt1}, \tiny{meg1}, \tiny{gemat11}, \tiny{bayer03}, \tiny{goodwin}},
ylabel={\tiny{CPU time (seconds)}},
ymajorgrids,
]
%% 1- SVD-cpu
\addplot+[smooth,color=black,loosely dotted, every mark/.append style={solid}, mark=+]
table[row sep=crcr]{
1	1.20000000000000  \\
2	6.59000000000000 \\
3	21.8300000000000 \\
4	54.9000000000000 \\
5	74.9700000000000  \\
};

%% 2- CPQR-cpu
\addplot+[smooth,color=gray,loosely dotted, every mark/.append style={solid}, mark=x]
table[row sep=crcr]{
1	0.77000000000000 \\
2	1.63000000000000 \\
3	8.36000000000000 \\
4	22.7700000000000 \\
5	29.5400000000000 \\
};

%% 3- QDWH-SVS-cpu
\addplot+[smooth,color=red, loosely dotted, every mark/.append style={solid}, mark=|]
table[row sep=crcr]{
1	6.46000000000000 \\ 
2	22.2300000000000 \\
3	70.0200000000000 \\
4	180.210000000000 \\
5	184.690000000000 \\
};

%% 4- randUTVq=0-cpu
\addplot+[smooth,color=blue,loosely dotted, every mark/.append style={solid}, mark=star]
table[row sep=crcr]{
1	7.07000000000000 \\
2	37.6400000000000 \\
3	191.140000000000 \\
4	470.480000000000 \\
5	591.010000000000 \\
};

%% 5- randUTVq=2-cpu
\addplot+[smooth,color=blue, loosely dotted, every mark/.append style={solid}, mark=asterisk]
table[row sep=crcr]{
1	7.33000000000000 \\
2	38.0400000000000 \\
3	193.590000000000 \\
4	481.400000000000 \\
5	592.950000000000  \\
};

%% 6- Rand-QLP-cpu
\addplot+[smooth,color=teal,loosely dotted, every mark/.append style={solid}, mark=-]
table[row sep=crcr]{
1	0.570000000000000 \\
2	1.38000000000000 \\
3	4.93000000000000 \\
4	10.1100000000000 \\
5	11.8800000000000  \\
};

\end{axis}

\begin{axis}[%
name=SumRate,
at={($(ber.east)+(35,0em)$)},
		anchor= west,
ymode=log,
width  = 0.32\columnwidth,%5.63489583333333in,
height = 0.24\columnwidth,%4.16838541666667in,
scale only axis,
xmin   = 1,
xmax  = 5,
xmajorgrids,
ymin = 0.28,
ymax = 577,
xtick       ={1, 2, 3, 4, 5},
xticklabels ={\tiny{watt1}, \tiny{meg1}, \tiny{gemat11}, \tiny{bayer03}, \tiny{goodwin}},
ymajorgrids,
ylabel={\tiny{Hybrid GPU time (seconds)}},
%title = {\texttt{HighNoiseLowRank}}]
%ytick       ={0.0973831, 0.0973830 , 0.0973829, 0.0973828},
%yticklabels ={$9.73831$, $9.73830$ , $9.73829$, $9.73828$},
legend entries={SVD,CPQR, QDWH-SVD, randUTV(q=0), randUTV(q=2), Rand-QLP}, legend style={at={(.6,1.2)},anchor=north east,draw=black,fill=white,legend cell align=left,font=\tiny, legend columns=3}
]

%% 1- SVD-gpu 
\addplot+[smooth,color=black,loosely dotted, every mark/.append style={solid}, mark=+]
table[row sep=crcr]{
1	0.800000000000000 \\
2	1.72000000000000 \\
3	7.33000000000000 \\
4	12.5600000000000 \\
5	18.1900000000000 \\
};

%% 2- CPQR-gpu
\addplot+[smooth,color=gray,loosely dotted, every mark/.append style={solid}, mark=x]
table[row sep=crcr]{
1	0.28000000000000 \\
2	1.40000000000000 \\
3	5.11000000000000 \\
4	10.1800000000000 \\
5	11.7900000000000  \\
};

%% 3- QDWH-SVD-gpu
\addplot+[smooth,color=red, loosely dotted, every mark/.append style={solid}, mark=|]
table[row sep=crcr]{
1	7.38000000000000 \\
2	19.3300000000000 \\
3	47.6400000000000 \\
4	132.100000000000 \\
5	140.140000000000 \\
};

%% 4- randUTVq=0-gpu
\addplot+[smooth,color=blue, loosely dotted, every mark/.append style={solid}, mark=star]
table[row sep=crcr]{
1	6.29000000000000 \\
2	34.9400000000000 \\
3	183.380000000000 \\
4	461.990000000000 \\
5	566.980000000000 \\
};

%% 5-randUTVq=2-gpu
\addplot+[smooth,color=blue, loosely dotted, every mark/.append style={solid}, mark=asterisk]
table[row sep=crcr]{
1	6.84000000000000 \\
2	37.5200000000000 \\
3	188.290000000000 \\ 
4	473.350000000000 \\
5	576.130000000000  \\
};

%% 6- Rand-QLP-gpu
\addplot+[smooth,color=teal,loosely dotted, every mark/.append style={thick}, mark=-]
table[row sep=crcr]{
1	0.48000000000000 \\
2	0.61000000000000 \\
3	1.73000000000000 \\
4	3.62000000000000 \\
5	4.41000000000000 \\
};

\end{axis}
\end{tikzpicture}%

\begin{tikzpicture}[font=\footnotesize] 

\begin{axis}[%
name=ber,
ymode=log,
width  = 0.32\columnwidth, 
height = 0.25\columnwidth, 
scale only axis,
xmin  = 1,
xmax  = 5,
xmajorgrids,
ymin = 0,
ymax = 51,
xtick       ={1, 2, 3, 4, 5},
xticklabels ={\tiny{watt1}, \tiny{meg1}, \tiny{gemat11}, \tiny{bayer03}, \tiny{goodwin}},
ylabel={\tiny{CPU speedup}},
ymajorgrids,
]
%% 1- /SVD-cpu
\addplot[smooth,mark=square*, black]
table[row sep=crcr]{
1	2.10530000000000 \\
2	4.77540000000000 \\
3	4.42800000000000 \\
4	5.43030000000000 \\
5	6.31060000000000 \\
};

%% 2- /CPQR-cpu
\addplot[smooth,mark=triangle*, gray]
table[row sep=crcr]{
1	1.35090000000000 \\
2	1.18120000000000 \\
3	1.69570000000000 \\
4	2.25220000000000 \\
5	2.48650000000000 \\
};

%% 3-/QDWH-SVD-cpu
\addplot[smooth,mark=*, red]
table[row sep=crcr]{
1	11.3333000000000 \\
2	16.1087000000000 \\
3	14.2028000000000 \\
4	17.8249000000000 \\
5	15.5463000000000 \\
};

%% 4- /randUTVq=0-cpu
\addplot[smooth,mark=diamond, blue]
table[row sep=crcr]{
1	12.4035000000000 \\
2	27.2754000000000 \\
3	38.7708000000000 \\
4	46.5361000000000 \\
5	49.7483000000000 \\
};

%% 5- /randUTVq=2-cpu
\addplot[smooth,mark=diamond*, blue]
table[row sep=crcr]{
1	12.8596000000000 \\
2	27.5652000000000 \\
3	39.2677000000000 \\
4	47.6162000000000 \\
5	49.9116000000000 \\
};

\end{axis}

\begin{axis}[%
name=SumRate,
at={($(ber.east)+(35,0em)$)},
		anchor= west,
ymode=log,
width  = 0.32\columnwidth, 
height = 0.25\columnwidth, 
scale only axis,
xmin   = 1,
xmax  = 5,
xmajorgrids,
ymin = 0.5,
ymax = 132,
xtick       ={1, 2, 3, 4, 5},
xticklabels ={\tiny{watt1}, \tiny{meg1}, \tiny{gemat11}, \tiny{bayer03}, \tiny{goodwin}},
ymajorgrids,
ylabel={\tiny{Hybrid GPU speedup}},
%title = {\texttt{HighNoiseLowRank}}]
%ytick       ={0.0973831, 0.0973830 , 0.0973829, 0.0973828},
%yticklabels ={$9.73831$, $9.73830$ , $9.73829$, $9.73828$},
legend entries={$t_\text{svd}/t_\text{randqlp}$,$t_\text{cpqr}/t_\text{randqlp}$, $t_\text{qdwhsvd}/t_\text{randqlp}$, $t_\text{randUTVq=0}/t_\text{randqlp}$, $t_\text{randUTVq=2}/t_\text{randqlp}$}, 
legend style={at={(0.8,1.20)},anchor=north east,draw=black,fill=white,legend cell align=left,font=\tiny, legend columns=3}
]

%% 1- /SVD gpu
\addplot[smooth,mark=square*, black]
table[row sep=crcr]{
1	1.66670000000000 \\
2	2.81970000000000 \\
3	4.23700000000000 \\
4	3.46960000000000 \\
5	4.12470000000000 \\
};

%% 2- /CPQR gpu
\addplot[smooth,mark=triangle*, gray]
table[row sep=crcr]{
1	0.58330000000000 \\
2	2.29510000000000 \\
3	2.95380000000000 \\
4	2.81220000000000 \\
5	2.67350000000000  \\
};

%% 3- /QDWH-SVD gpu
\addplot[smooth,mark=*, red]
table[row sep=crcr]{
1	15.3750000000000 \\
2	31.6885000000000 \\
3	27.5376000000000 \\
4	36.4917000000000 \\
5	31.7778000000000 \\
};

%% 4- /randUTVq=0 gpu
\addplot[smooth,mark=diamond, blue]
table[row sep=crcr]{
1	13.1042000000000 \\
2	57.2787000000000 \\
3	106 \\
4	127.621500000000 \\
5	128.566900000000 \\
};

%%%%%%%%%%%
%% 5- /randUTVq=2 gpu
\addplot[smooth,mark=diamond*, blue]
table[row sep=crcr]{
1	14.2500000000000 \\
2	61.5082000000000 \\
3	108.838200000000 \\
4	130.759700000000 \\
5	130.641700000000 \\
};

\end{axis}
\end{tikzpicture}%     
\captionsetup{justification=centering,font=scriptsize}  
\caption{Computational times for different algorithms on real matrices (top row) and speedups offered by Rand-QLP (bottom row). For randUTV, we set the block size $b$=100 and oversampling parameter to 5.} 
\label{figSpReal}      
\end{center}
\end{figure}
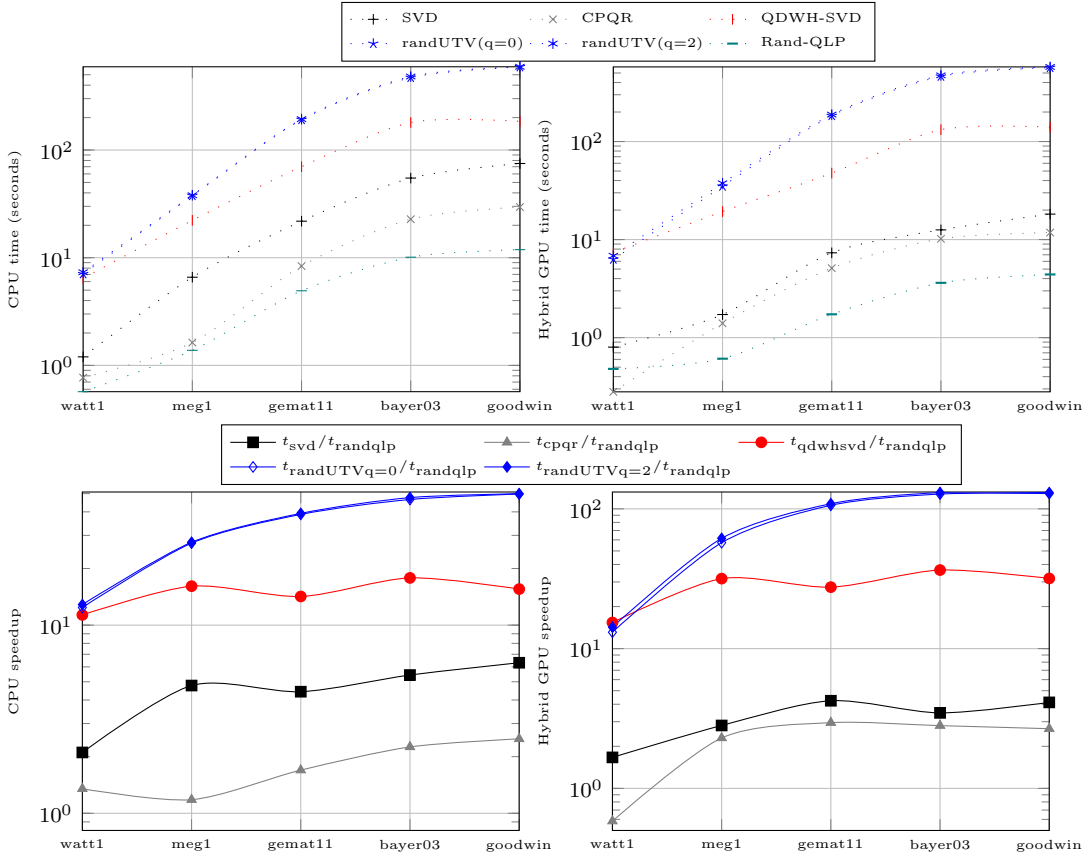

Computational times for different algorithms and speedups offered by Rand-QLP are shown in Figures \ref{figSpRandom} and \ref{figSpReal}. The results for randUTV and Rand-QLP were averaged over 5 runs, in which we observe that Rand-QLP is substantially outperforms the other algorithms. Specifically: 
\begin{itemize}
\item Uing only the CPU, for random matrices, Rand-QLP achieves speedups of up to 6.6$\times$, 5$\times$, 23.6$\times$, and 60.8$\times$ over the SVD, CPQR, QDWH-SVD, and randUTV. For real matrices, these speedups are 6.3$\times$, 2.4$\times$, 17.8$\times$, and 49.7$\times$, repsetively.
\item Employing the GPU, for random matrices, Rand-QLP achieves speedups of up to 4.4$\times$, 3.8$\times$, 34.2$\times$, and 154.4$\times$ over the SVD, CPQR, QDWH-SVD, and randUTV. These speedups for real matrices are 4.2$\times$, 2.8$\times$, 36.2$\times$, and 130.6$\times$, respectively. 
\end{itemize}

\subsection{Singular Values and Low-Rank Approximations}

We compare the performance of the algorithms considered for the tasks of singular values and low-rank approximations on randomly generated matrices, formed as described below, as well as the real matrices presented in Table \ref{TableMatSuiSpa}. 

We generate four types of low-rank matrices of order $n=1000$. The first matrix, named \texttt{LowRankPlusNoise}, takes the following form \cite{DemGGX15}: 
 \begin{equation}\label{eq_A1PlusA2}
  {\bf A} = {\bf A}_\text{LR} + {\bf A}_n.
 \end{equation}
Here ${\bf A}_\text{LR} = {\bf U\Sigma V}^T$, where ${\bf U}, {\bf V} \in \mathbb R^{n \times n}$ are random orthogonal matrices, and ${\bf \Sigma} \in \mathbb R^{n \times n}$ is diagonal with entries $\sigma_i$s as singular values. The singular values decrease linearly from $1$ to $10^{-20}$, and we set $\sigma_{k+1}=...=\sigma_n=0$. The matrix ${\bf A}_n=0.05\sigma_k{\bf G}_N$, where ${\bf G}_N$ is normalized Gaussian. We set $k=200$.  

The other three matrices are formed as ${\bf A}_\text{LR}$, but the distribution of their singular values take the following form: 
\begin{itemize}
\item \texttt{LowRankFastDecay} \cite{YuLi18}: ${\bf \Sigma}= \text{diag}(i^{-2})$ for $i=1, ..., n$.
\item \texttt{LowRankS-ShapedDecay} \cite{MartinssonQH19}: The entries of $\bf \Sigma$ hover around 1, then decay quickly and finally level out at $0.01$.
\item \texttt{LowRankSlowDecay} \cite{TrYUC17}: ${\bf \Sigma}= \text{diag}(\underbrace{
1, ..., 1}_{k=100}, 2^{-1}, 3^{-1}, ..., (n-k+1)^{-1})$. 
\end{itemize}  

\subsubsection{Singular values approximation} 
The results of approximating the singular values of $\bf A$ are shown in Figures \ref{figSVsyn} and \ref{figSSparse}. We observe in Figure \ref{figSVsyn} that for \texttt{LowRankPlusNoise}, Rand-QLP strongly reveals the gap in the spectrum and provides approximations that are very close to those of p-QLP and the optimal SVD, while CPQR only suggests the gap and underestimates the leading singular values. randUTV with no power iteration ($q=0$) weakly reveals the rank and also introduces extra gaps in its approximations, thereby underestimating or overestimating most singular values. This result was expected: randUTV treats each block (of b columns) of an input matrix separately, and randomized algorithms tend to overestimate or/and underestimate singular values of a matrix if the power iteration is not used (see, e.g., \cite{MFKDeTSP18}). In spite of improvement in approximation accuracy of randUTV by setting $q=2$, the result suggests that a larger value of $q$ is required for obtaining highly accurate approximations to the singular values. This in turn increases the cost of randUTV computation.  

For \texttt{LowRankFastDecay}, \texttt{LowRankS-ShapedDecay}, and \texttt{LowRankSlowDecay}, CPQR tends to underestimate or/and overestimate the singular values, randUTV shows a similar behavior as for \texttt{LowRankPlusNoise}, and the estimations by Rand-QLP are indistinguishable from those computed by p-QLP and the SVD. 

In Figure \ref{figSSparse}, we observe that CPQR provides more accurate estimations to the singular values of the real matrices, while the estimations by Rand-QLP are very close to those of the SVD. 

\subsubsection{Low-rank approximation}

For each matrix $\bf A$, we compute rank-$k$ approximations for a general $1 \le k \le n$ as follows:
\begin{equation}\notag
\begin{aligned}
& {\bf A}_{k} = {\bf U}(:, 1:k){\bf U}(1:k, 1:k){\bf V}(:, 1:k)^T. \hspace{.2cm} \text{\% SVD} \\
& {\bf A}_{k} = {\bf Q}_\text{piv}(:, 1:k){\bf R}_\text{piv}(1:k,:){\bf \Pi}^T. \hspace{1.2cm} \text{\% CPQR} \\
& {\bf A}_{k} = {\bf Q}_\text{qlp}{\bf L}_\text{qlp}(:, 1:k){\bf P}_\text{qlp}(1:k, :)^T. \hspace{.9cm} \text{\% Pivoted QLP}\\
& {\bf A}_{k} = {\bf U}(:, 1:k){\bf T}(1:k,:){\bf V}^T. \hspace{2cm} \text{\% randUTV}\\
& {\bf A}_{k} = {\bf Q}{\bf L}(:, 1:k){\bf P}(1:k, :)^T. \hspace{2.1cm} \text{\% Rand-QLP}
\end{aligned}
\end{equation}

\begin{figure}[!tbp]
\captionsetup[subfigure]{labelformat=empty}
  \centering
  \subfloat[]{\includegraphics[width=0.49\textwidth]{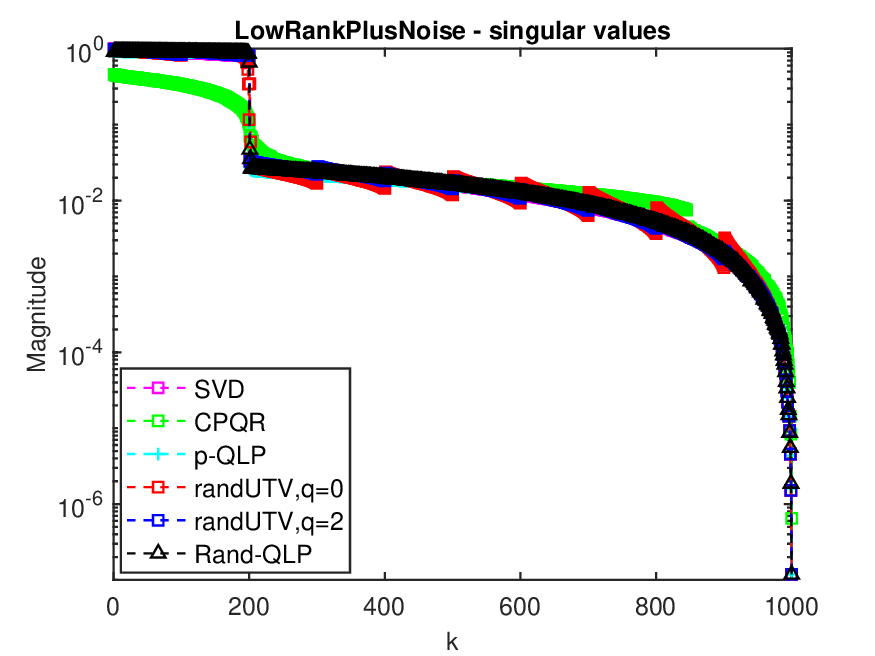}}
  \hfill
  \subfloat[]{\includegraphics[width=0.49\textwidth]{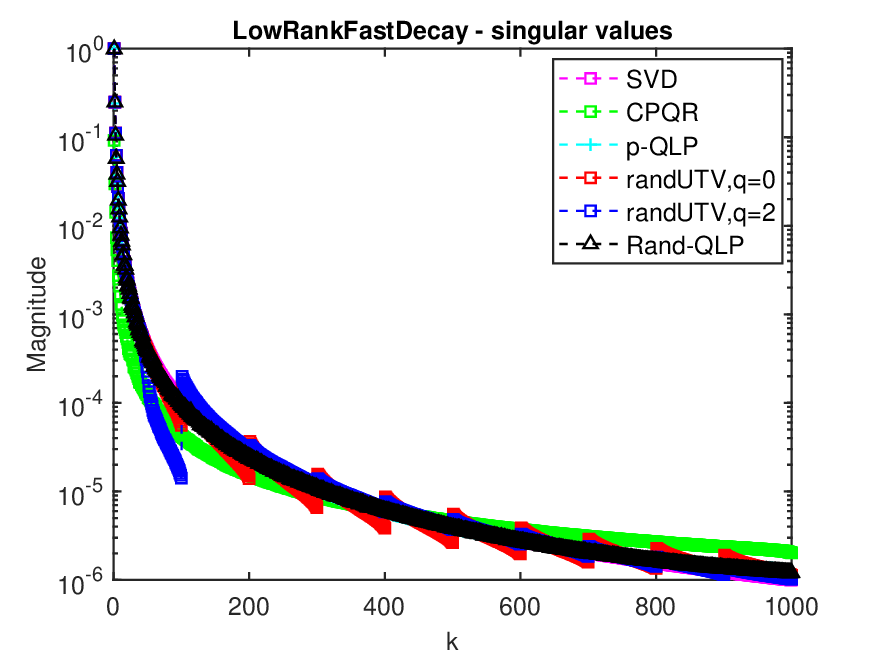}}
  \hfill
  \subfloat[]{\includegraphics[width=0.49\textwidth]{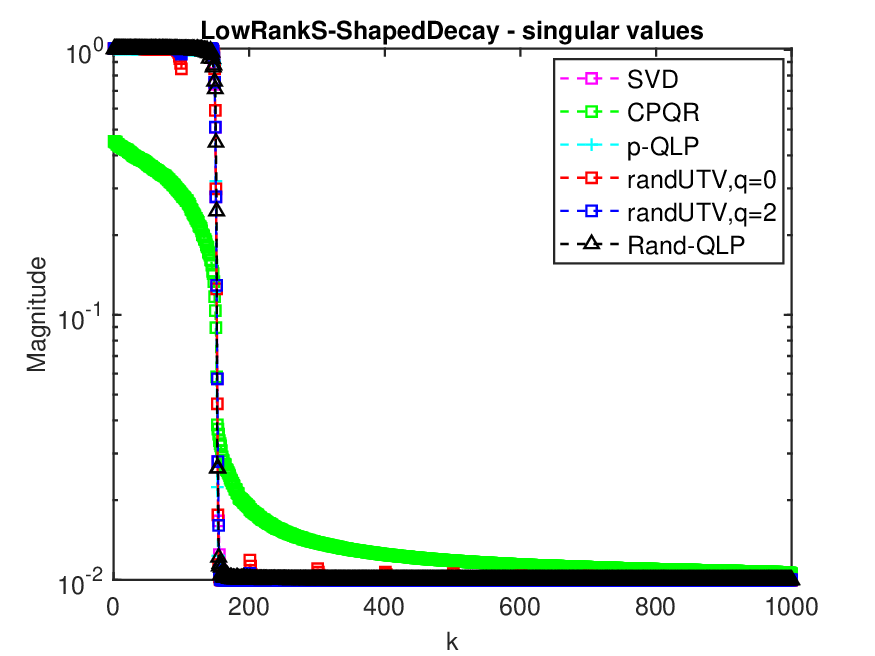}}
 \hfill
  \subfloat[]{\includegraphics[width=0.49\textwidth]{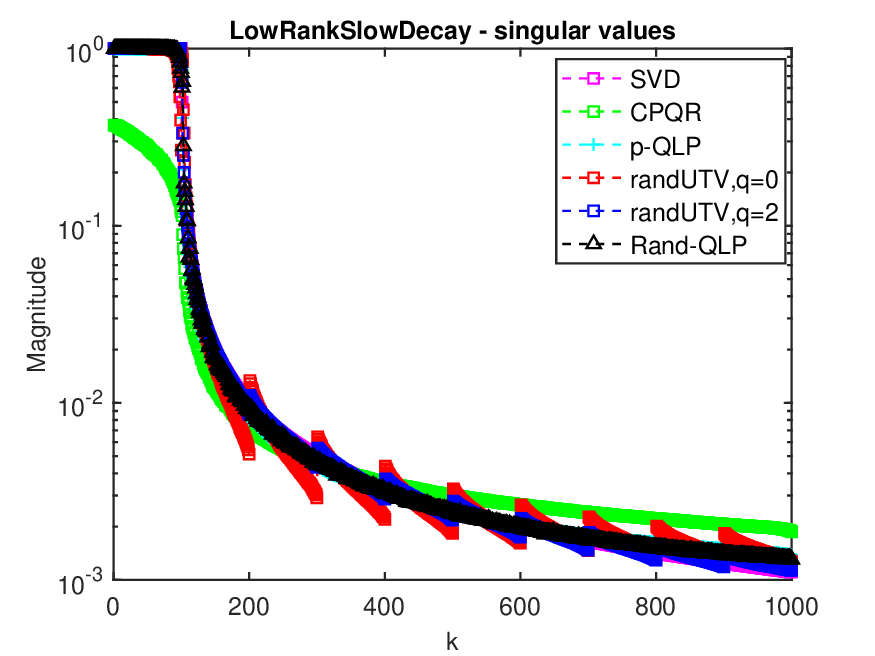}} 
  \vspace{-0.8cm}
  \caption{Estimates of the singular values for \texttt{LowRankPlusNoise}, \texttt{LowRankFastDecay}, \texttt{LowRankS-ShapedDecay}, and \texttt{LowRankSlowDecay}.}
  \label{figSVsyn}
\end{figure}

\begin{figure}[!tbp]
\captionsetup[subfigure]{labelformat=empty}
  \centering
  \subfloat[]{\includegraphics[width=0.49\textwidth]{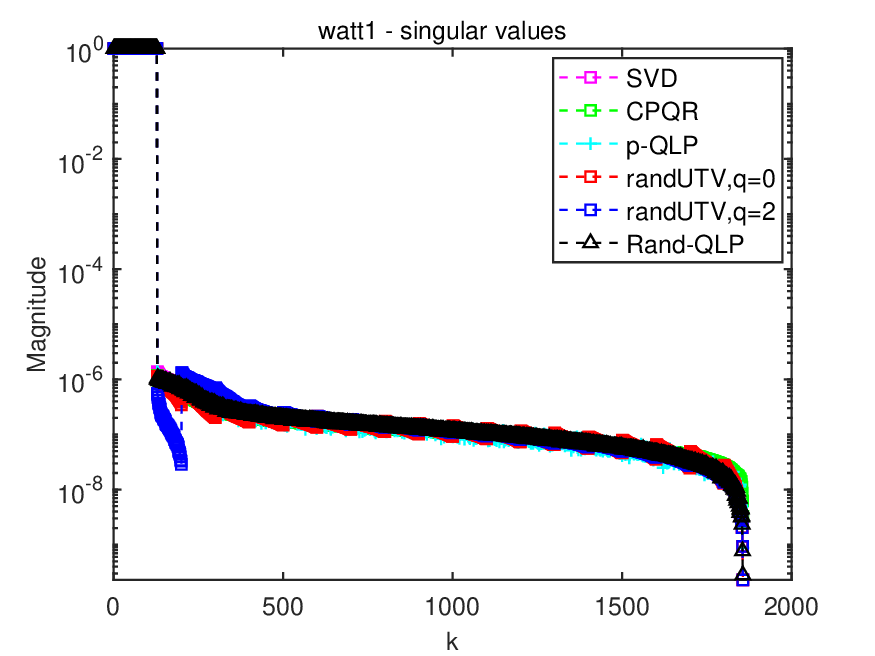}}
  \hfill
  \subfloat[]{\includegraphics[width=0.49\textwidth]{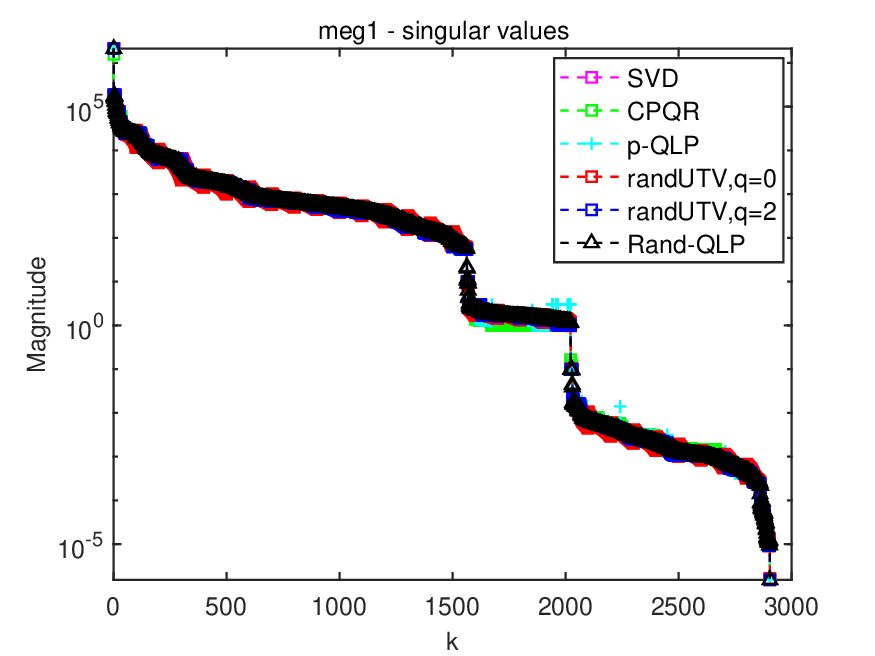}}
  \hfill
  \subfloat[]{\includegraphics[width=0.49\textwidth]{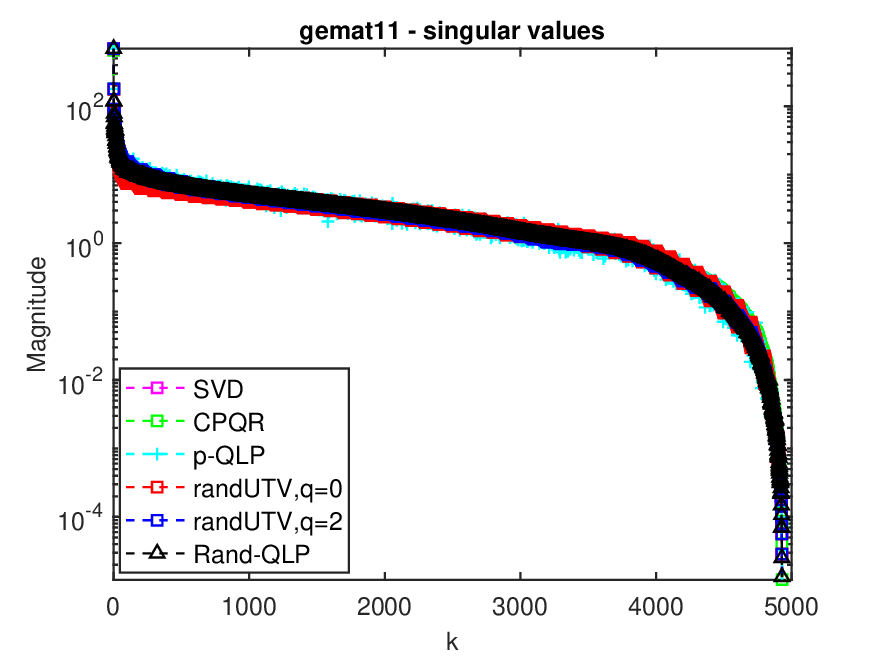}}
  \hfill
  \subfloat[]{\includegraphics[width=0.49\textwidth]{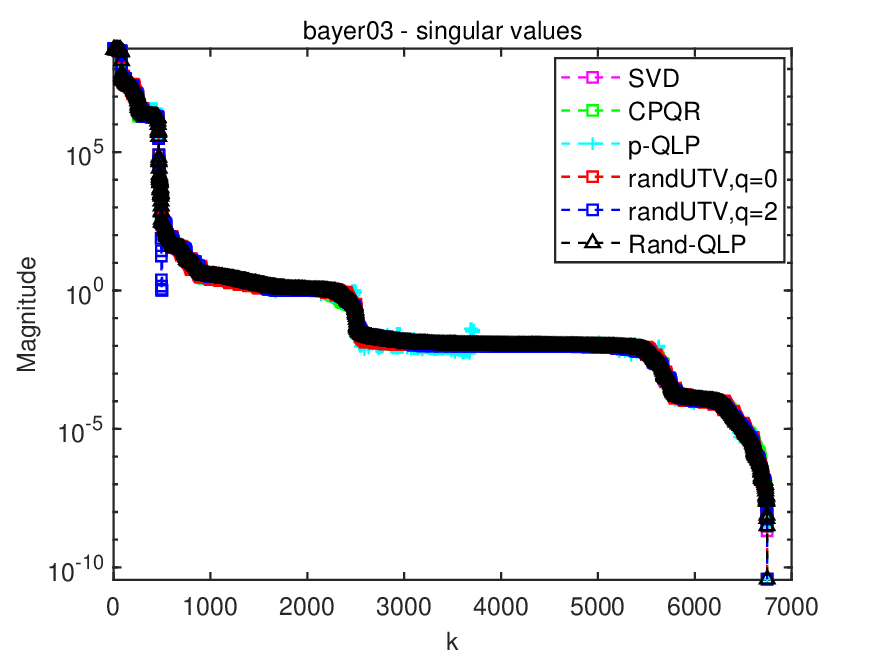}}
  \hfill
  \subfloat[]{\includegraphics[width=0.49\textwidth]{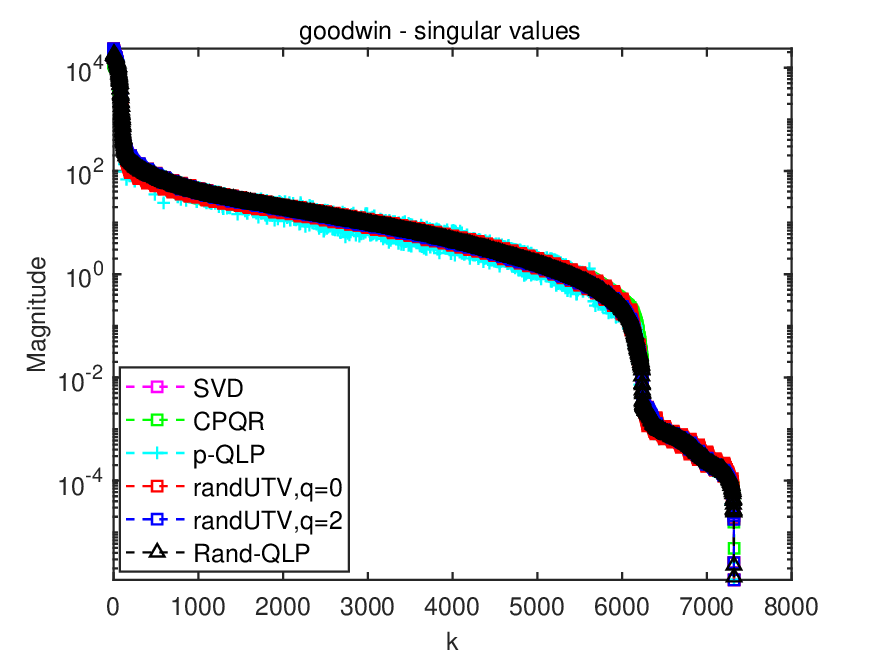}}
  \vspace{-0.8cm}
  \caption{Estimates of the singular values for real matrices from the SuiteSparse Matrix Collection presented in Table \ref{TableMatSuiSpa}.}
  \label{figSSparse}
\end{figure}

\begin{figure}[!tbp]
\captionsetup[subfigure]{labelformat=empty}
  \centering
  \subfloat[]{\includegraphics[width=0.4\textwidth]{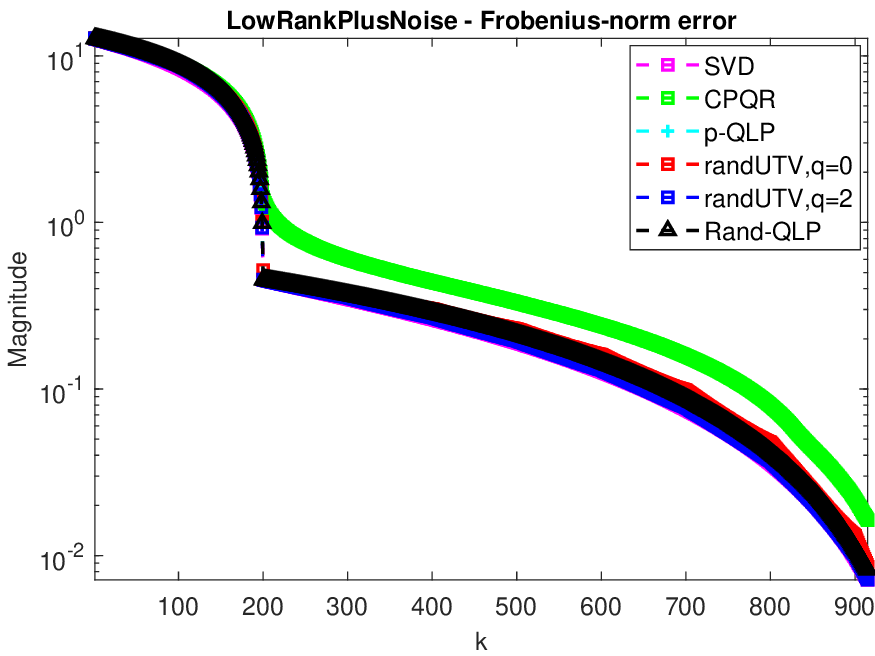}}
  \hfill
\subfloat[]{\includegraphics[width=0.4\textwidth]{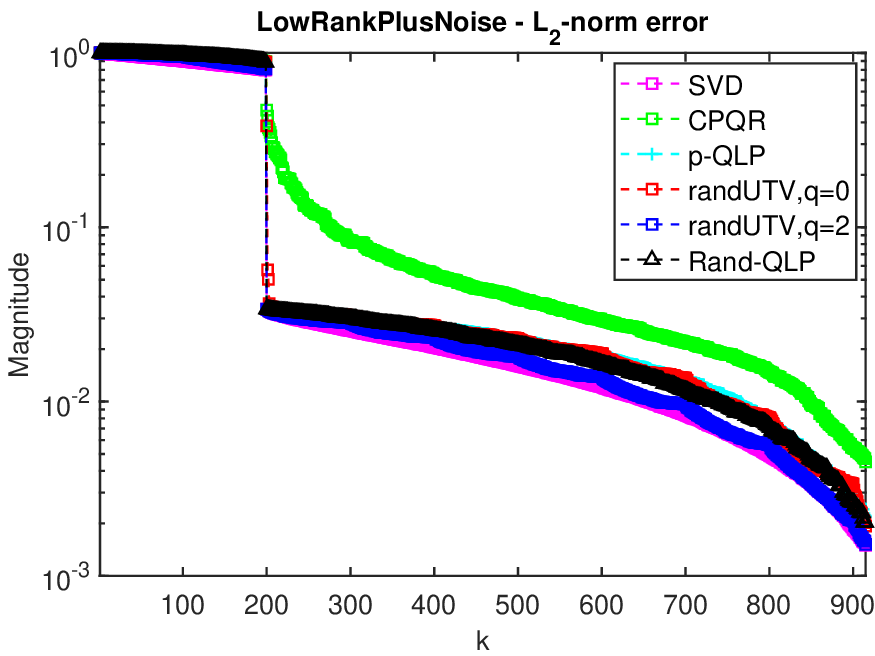}}
  \hfill
  \subfloat[]{\includegraphics[width=0.4\textwidth]{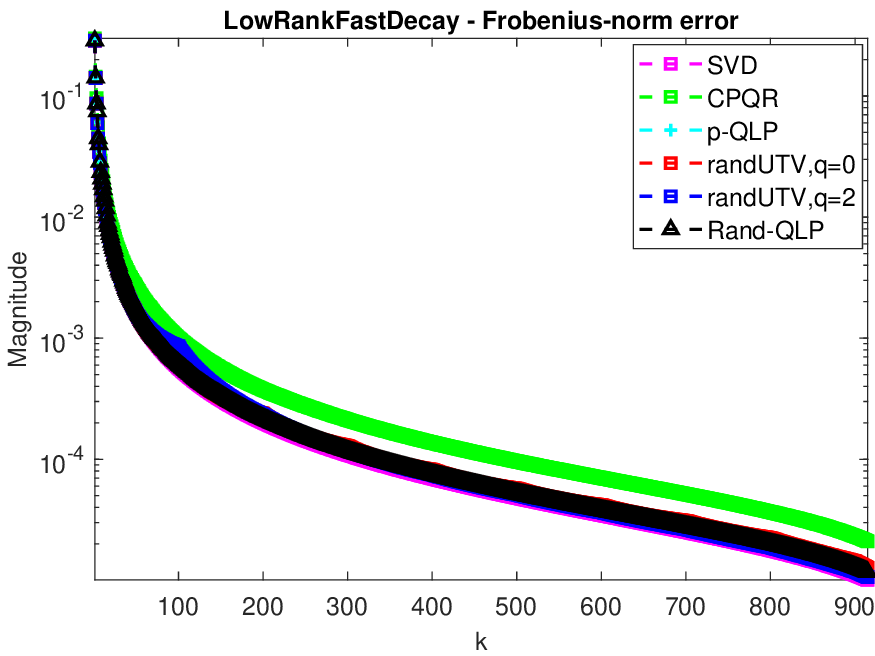}}
  \hfill
  \subfloat[]{\includegraphics[width=0.4\textwidth]{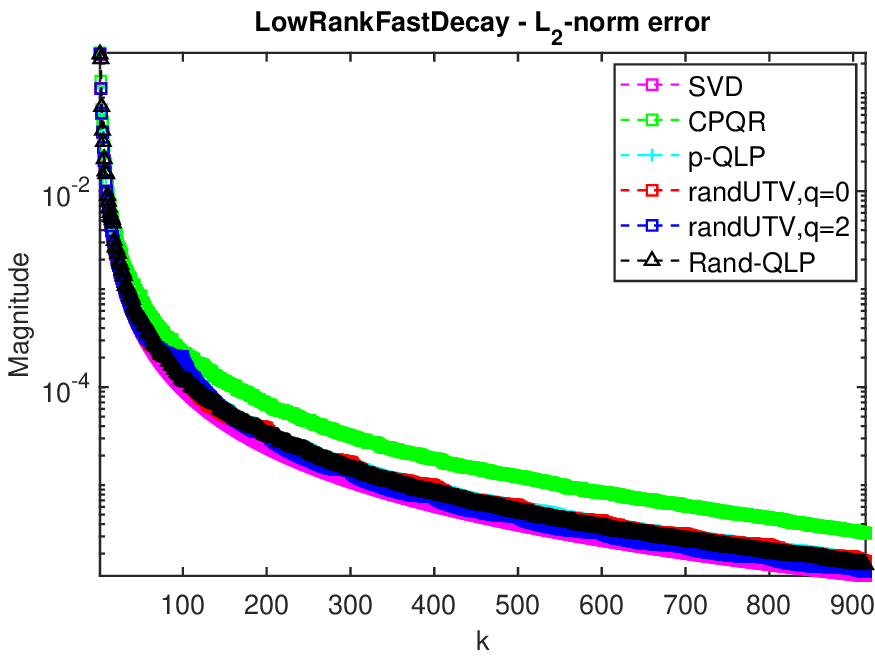}}
  \hfill
  \subfloat[]{\includegraphics[width=0.4\textwidth]{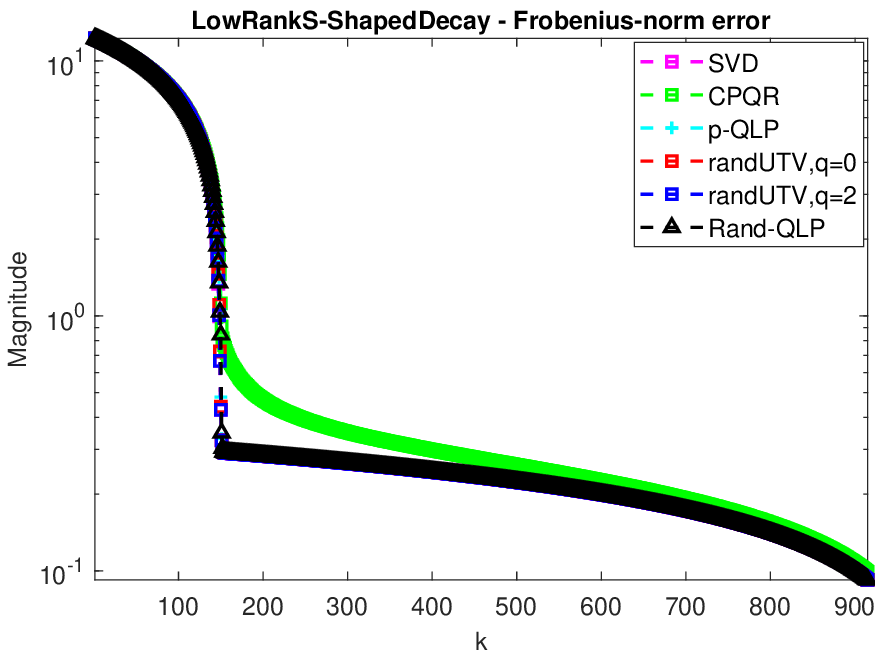}}
  \hfill
  \subfloat[]{\includegraphics[width=0.4\textwidth]{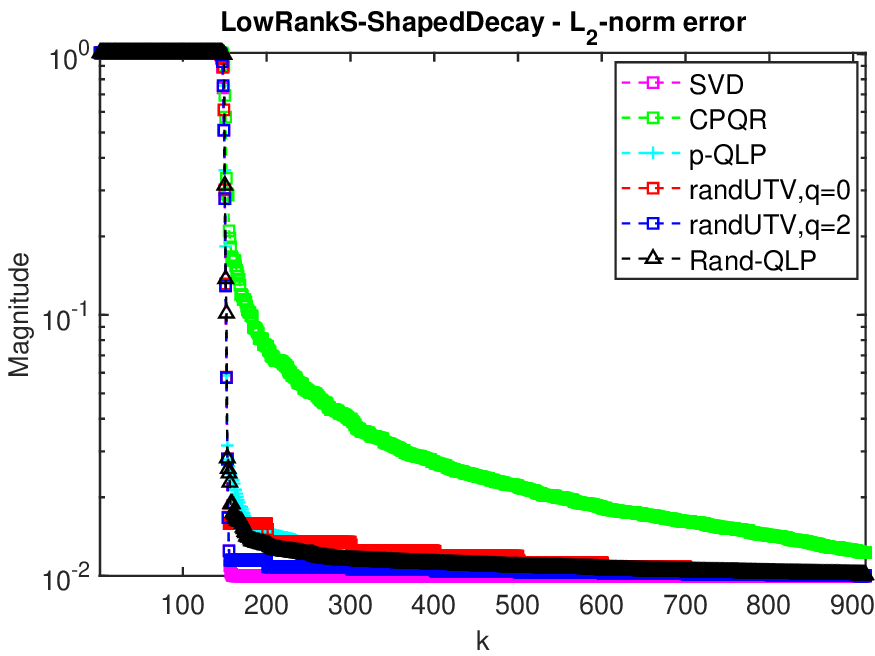}}
\hfill
  \subfloat[]{\includegraphics[width=0.4\textwidth]{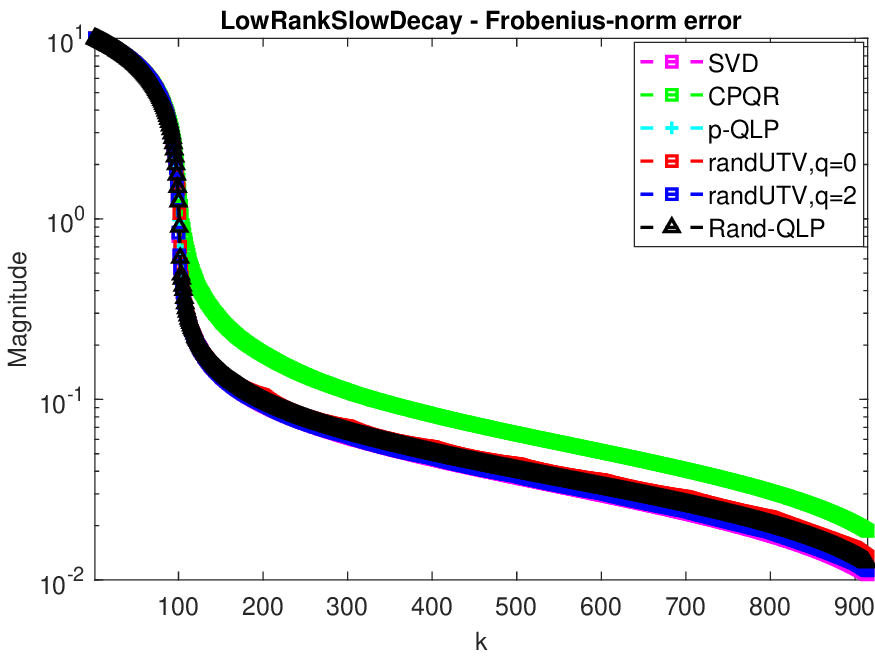}}
  \hfill
  \subfloat[]{\includegraphics[width=0.4\textwidth]{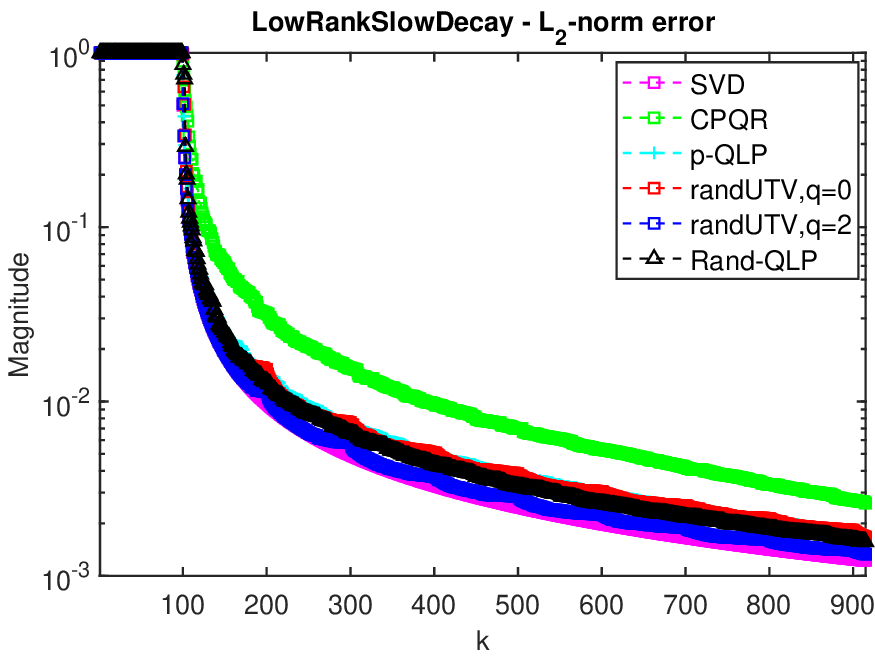}}
  \vspace{-0.8cm}
  \caption{Frobenius and $\ell_2$-norm approximation errors for \texttt{LowRankPlusNoise}, \texttt{LowRankFastDecay}, \texttt{LowRankS-ShapedDecay}, and \texttt{LowRankSlowDecay}.}
  \label{figLRAppSy}
\end{figure}

\begin{figure}[!tbp]
\captionsetup[subfigure]{labelformat=empty}
  \centering
 \subfloat[]{\includegraphics[width=0.33\textwidth]{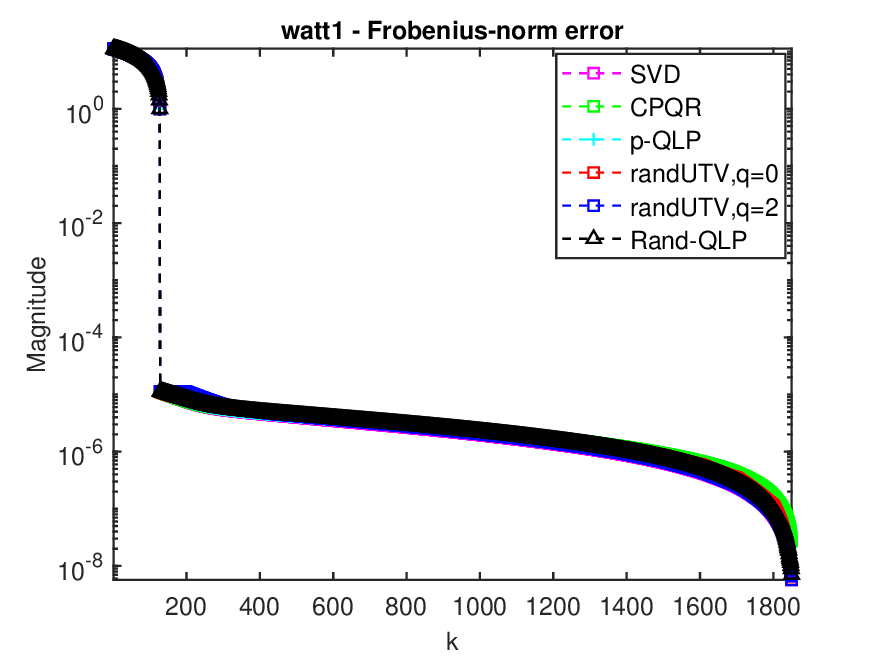}}
  %\hfill
\subfloat[]{\includegraphics[width=0.33\textwidth]{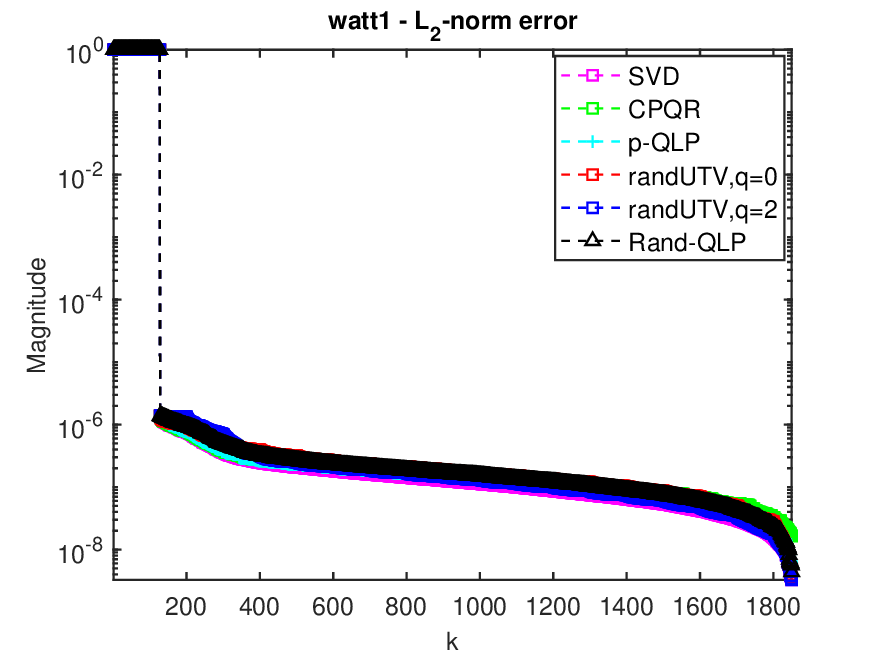}}
  %\hfill
  \subfloat[]{\includegraphics[width=0.33\textwidth]{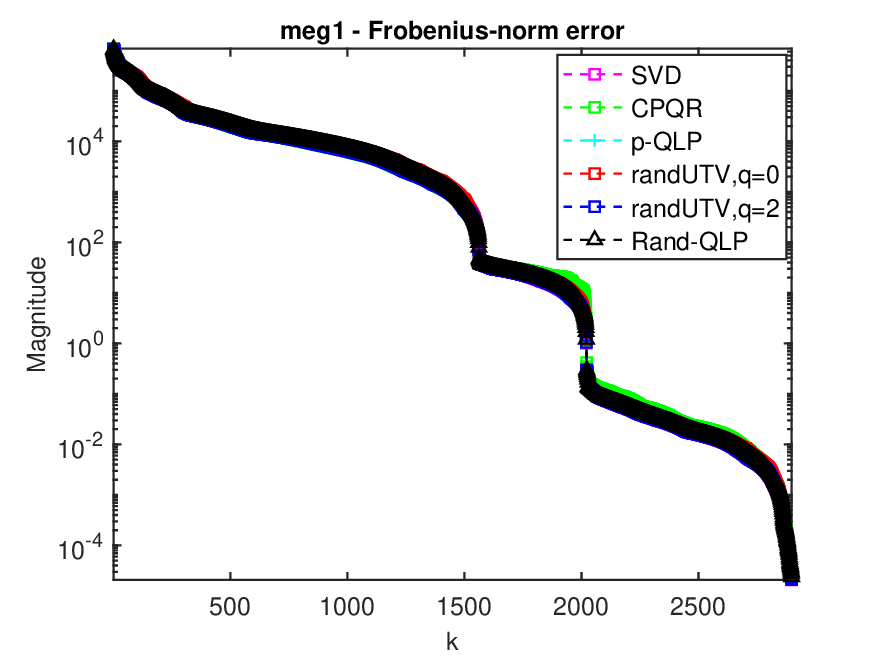}}
  \hfill
  \subfloat[]{\includegraphics[width=0.33\textwidth]{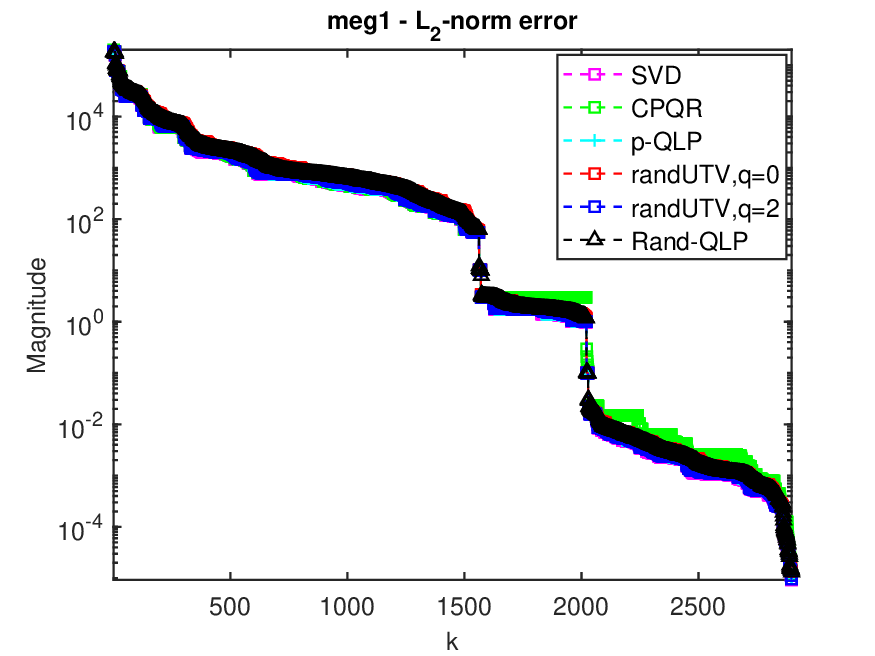}}
  %\hfill
  \subfloat[]{\includegraphics[width=0.33\textwidth]{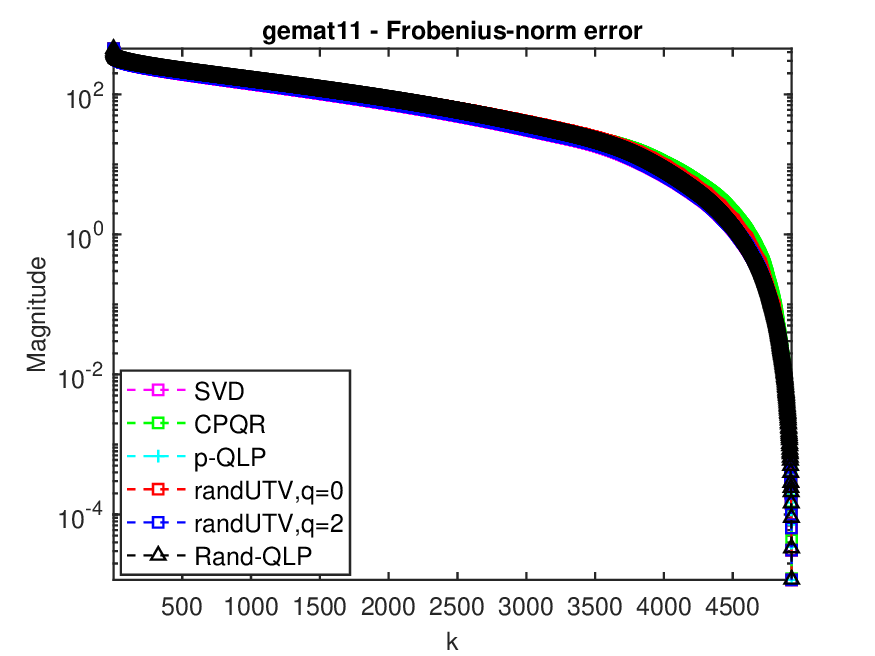}}
  %\hfill
  \subfloat[]{\includegraphics[width=0.33\textwidth]{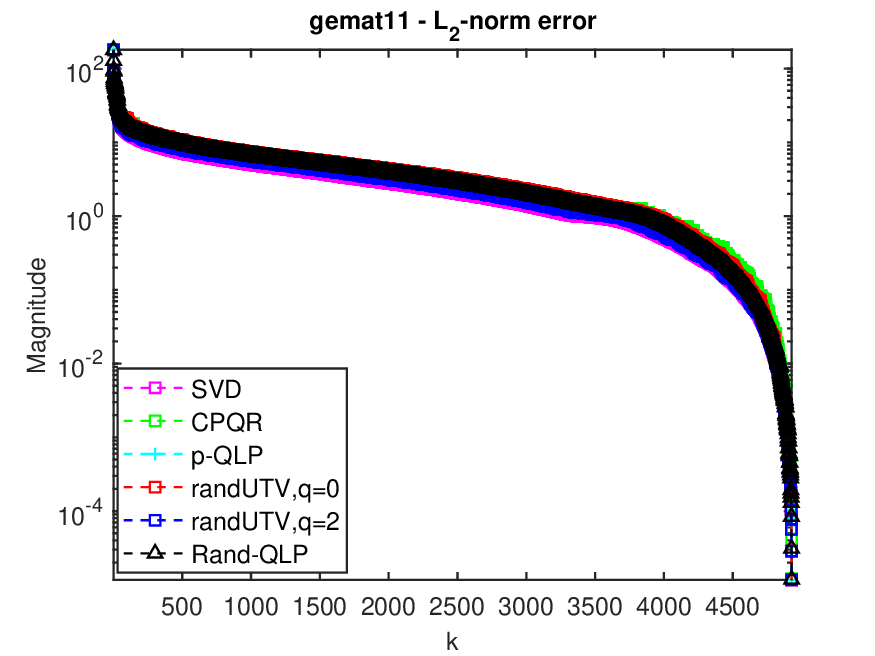}}
  \hfill
  \subfloat[]{\includegraphics[width=0.33\textwidth]{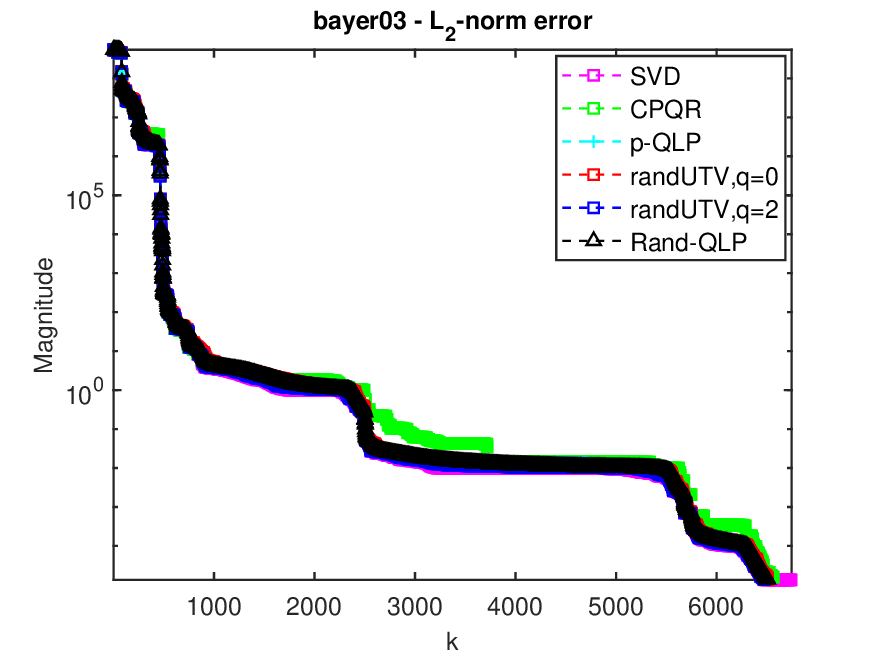}}
  %\hfill
  \subfloat[]{\includegraphics[width=0.33\textwidth]{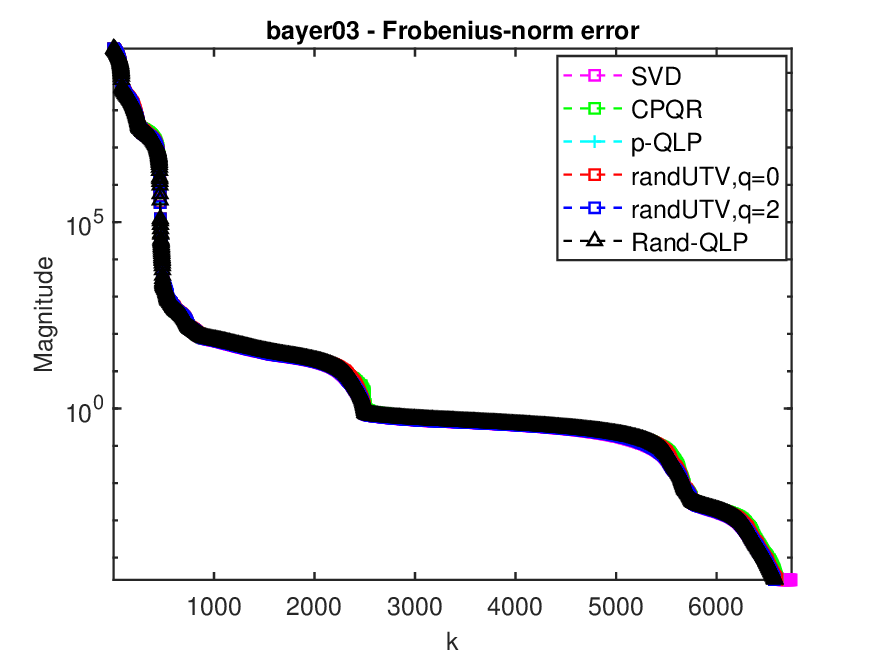}}
  %\hfill
  \subfloat[]{\includegraphics[width=0.33\textwidth]{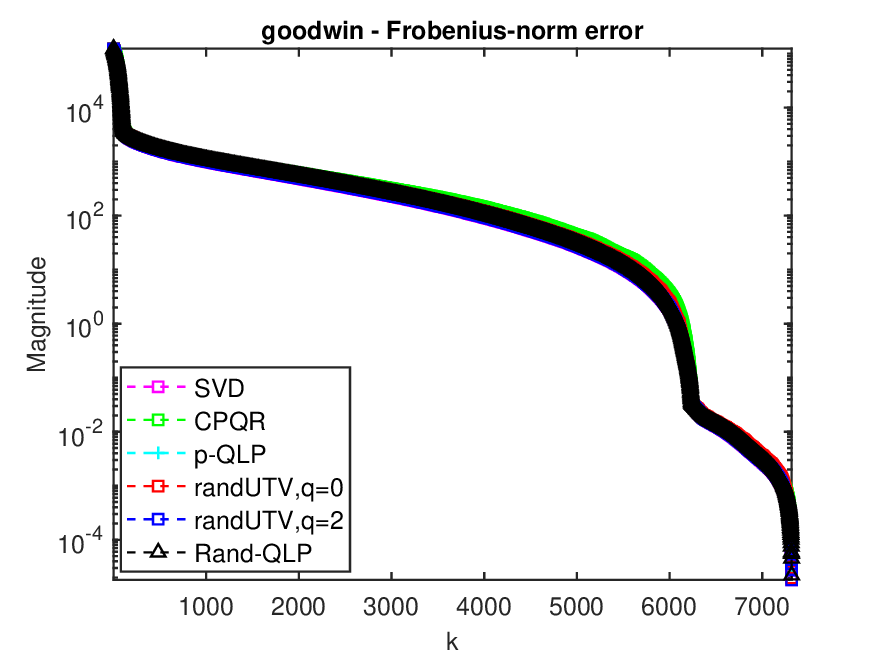}}
  \hfill
  \subfloat[]{\includegraphics[width=0.33\textwidth]{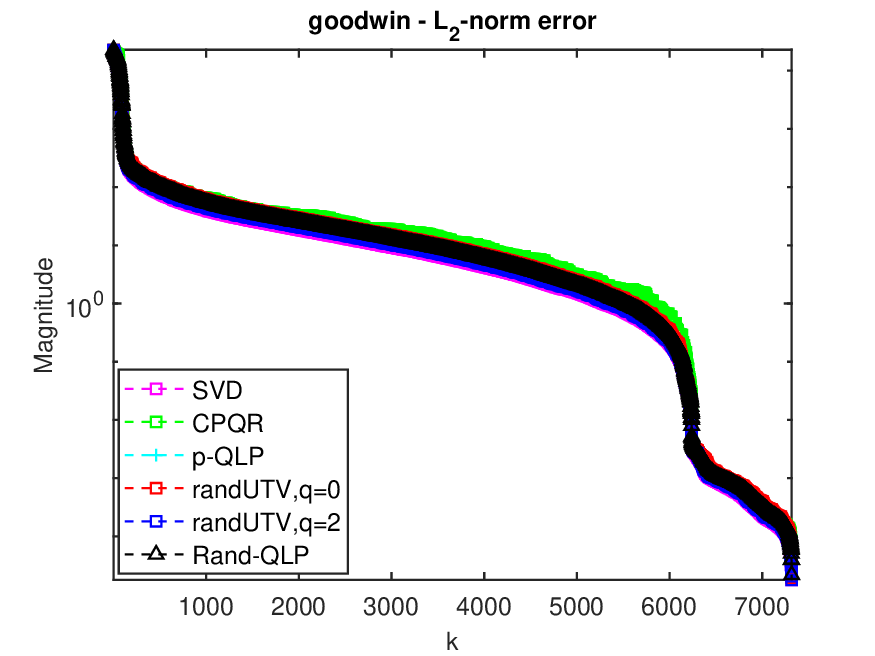}}
  \vspace{-0.8cm}
  \caption{Frobenius and $\ell_2$-norm approximation errors for real matrices from the SuiteSparse Matrix Collection presented in Table \ref{TableMatSuiSpa}.}
  \label{figLRAppSSparse}
\end{figure}

We then compute the Frobenius error $\|{\bf A} - {\bf A}_k\|_F$ as well as $\ell_2$-norm error $\|{\bf A} - {\bf A}_k\|_2$. The results are plotted in Figures \ref{figLRAppSy} and \ref{figLRAppSSparse}. On synthetic matrices, CPQR shows the poorest performance, while approximation errors produced by Rand-QLP are comparable to those incurred by p-QLP and the SVD. On real matrices, although all algorithms show similar behavior, Rand-QLP is computationally the most efficient algorithm and lends itself to parallelization much better than other algorithms, as demonstrated in Figures \ref{figSVsyn} and \ref{figSSparse}.

%%%%%%%%%%%%%%% singular values - synthetic data

%%%%%%%%%%%%%%% singular values - real data

%%%%%%%%%%%%%%% Synthetic Rank-\ell approximation

%%%%%%%%%%%%%%% REAL Rank-\ell approximation

\vspace{1.3cm}
\section{Conclusion}
\label{secCon}
In this paper we presented Rand-QLP, a randomized rank-revealing algorithm that approximates the full SVD of a matrix. Its computation requires only matrix-matrix multiplication and the unpivoted QR factorization and can be done almost entirely using level-3 BLAS routines. Rand-QLP can therefore be efficiently parallelized on modern computational platforms. We derived error bounds on i) the approximate leading singular values, ii) the distance between approximate four fundamental subspaces and the exact ones, and iii) the low-rank approximations. We showed through numerical tests on a hybrid GPU-accelerated multicore system that Rand-QLP can be substantially faster than existing algorithms, while providing approximations that are very close to those of p-QLP and the optimal SVD. 

%%%%%%%%%%%%%%%%%%%%%%%%%%%%%%%%%%%%%%%%%%%%%%%%%%%%%%%%%%%%%%%%%%%%%%%%%%

\section{Appendix}

This appendix provides proofs for the results presented in Theorems \ref{ThSinRQLP}, \ref{ThSins} and \ref{ThLowRankErr}. 

\subsection{Proof of Theorem \ref{ThSinRQLP}}
\label{AppdxPrTh1}

\subsubsection{Proof of equation \eqref{eqSinV}}

Let the matrix $\bf A$ have a Rand-QLP decomposition as stated in equation (3). We have for the range of the matrix $\bf Q$ in Step 4 of Algorithm 1:
\begin{equation}\notag
\mathcal{R}({\bf Q}) = \mathcal{R}({\bf A}\bar{\bf Q})= \mathcal{R}({\bf AA}^T{\bf \Omega}).
\end{equation}

Let matrices ${\bf Q}$ and $\bar{\bf Q}$ (Step 3 of Algorithm 1) be partitioned as follows: %%%%%%%%%%%%%%
\begin{equation}\notag
{\bf Q} =[{\bf Q}_1\quad {\bf Q}_2], \quad \text{and} \quad \bar{\bf Q} =[\bar{\bf Q}_1\quad \bar{\bf Q}_2],
\end{equation}
where ${\bf Q}_1$ and $\bar{\bf Q}_1$ contain the first $k$ columns of the corresponding matrices. Hence 
\begin{equation}\notag
\mathcal{R}({\bf Q}_1) = \mathcal{R}({\bf A}\bar{\bf Q}_1)= \mathcal{R}({\bf AA}^T
{\bf \Omega}_1).
\end{equation}
Upon substitution of $\bf A$ by its SVD, we obtain %%%%%%%%%%%%%%
\begin{equation}\notag
\begin{aligned}
{\bf AA}^T{\bf \Omega}_1 & = {\bf U}
\begin{bmatrix}
{\bf \Sigma}_k^2\widetilde{\bf \Omega}_{1}\\
{\bf \Sigma}_\perp^2\widetilde{\bf \Omega}_{2}
\end{bmatrix}  = \ddot{\bf Q}\ddot{\bf R}.
\end{aligned}
\end{equation}
We thus have $\mathcal{R}({\bf Q}_1)=\mathcal{R}(\ddot{\bf Q})$, and both $\widetilde{\bf \Omega}_{1}$ and $\widetilde{\bf \Omega}_{2}$ have the standard Gaussian distribution due to the rotationally invariant property of $\bf \Omega$. We now define a $k\times k$ matrix ${\bf X} \triangleq \widetilde{\bf \Omega}_{1}^{-1}{\bf \Sigma}_k^{-2}$ and form the following matrix product on which a QR factorization is performed %%%%%%%%%%%%%%
\begin{equation} \label{eq_QRtilde}
\begin{aligned}
{\bf AA}^T{\bf \Omega}_1{\bf X} = {\bf U}\begin{bmatrix}
{\bf I}_k  \\
{\bf \Sigma}_\perp^2\widetilde{\bf \Omega}_{2}\widetilde{\bf \Omega}_{1}^{-1} {\bf \Sigma}_k^{-2}
\end{bmatrix}=\widetilde{\bf Q}\widetilde{\bf R}.
\end{aligned}
\end{equation}
Let ${\bf G} \triangleq {\bf \Sigma}_\perp^2\widetilde{\bf \Omega}_{2}\widetilde{\bf \Omega}_{1}^{-1} {\bf \Sigma}_k^{-2}$, and %%%%%%%%%%%%%%
\begin{equation}\notag
\begin{aligned}
{\bf F}^{-1} \triangleq \widetilde{\bf R}^{-1} \widetilde{\bf R}^{-T} = 
(\widetilde{\bf R}^T \widetilde{\bf R})^{-1} = ({\bf I}+{\bf G}^T{\bf G}) ^{-1}.
\end{aligned}
\end{equation}
Then, for $\widetilde{\bf Q}\widetilde{\bf Q}^T$ we obtain %%%%%%%%%%%%%%
\begin{equation}\notag
\widetilde{\bf Q}\widetilde{\bf Q}^T = {\bf U}
\begin{bmatrix}
{\bf F}^{-1} & {\bf F}^{-1}{\bf G}^T  \\
{\bf GF}^{-1} & {\bf G}{\bf F}^{-1}{\bf G}^T 
\end{bmatrix}{\bf U}^T.
\end{equation}
By algebraic equivalency, we have %%%%%%%%%%%%%%
\begin{equation}\notag
{\bf A}^T{\bf Q}_1  = {\bf P}_1{\bf L}_{11}^T,
\end{equation}
and for $i=1,..., k$, $\sigma_i({\bf A}^T{\bf Q}_1) = \sigma_i({\bf P}_1{\bf L}_{11}^T)$. 
From ${\bf I} \succeq {\bf Q}_1{\bf Q}_1^T$, we have ${\bf A}^T{\bf A} \succeq {\bf A}^T{\bf Q}_1{\bf Q}_1^T{\bf A}$, which, by leveraging the Cauchy's interlacing theorem \cite{StewartSun90}, gives %%%%%%%%%%%%%%
\begin{equation}\label{eqQ1eqQtilde}
\lambda_i({\bf A}^T{\bf A}) \ge \lambda_i({\bf A}^T{\bf Q}_1{\bf Q}_1^T{\bf A}) = \lambda_i({\bf A}^T\ddot{\bf Q}\ddot{\bf Q}^T{\bf A}) = \lambda_i({\bf A}^T\widetilde{\bf Q}\widetilde{\bf Q}^T{\bf A}).
\end{equation}
The last relation follows because for an orthogonal matrix $\bf W$ of appropriate size, we have $\ddot{\bf Q} = \widetilde{\bf Q}{\bf W}$ (the same principle expresses the relation between ${\bf Q}_1$ and $\ddot{\bf Q}$). Upon substitution we obtain %%%%%%%%%%%%%%%%%
\begin{equation}\notag
{\bf A}^T\widetilde{\bf Q}\widetilde{\bf Q}^T{\bf A} = {\bf V}
\begin{bmatrix}
{\bf \Sigma}_k{\bf F}^{-1}{\bf \Sigma}_k & {\bf \Sigma}_k{\bf F}^{-1}{\bf G}^T{\bf \Sigma}_\perp  \\
{\bf \Sigma}_\perp{\bf GF}^{-1}{\bf \Sigma}_k & {\bf \Sigma}_\perp{\bf G}{\bf F}^{-1}{\bf G}^T{\bf \Sigma}_k 
\end{bmatrix}{\bf V}^T.
\end{equation}
Therefore %%%%%%%%%%%%%%%%%%
\begin{equation}\notag
\lambda_i({\bf A}^T{\bf A}) \ge \lambda_i({\bf A}^T\widetilde{\bf Q}\widetilde{\bf Q}^T{\bf A}) \ge \lambda_i({\bf \Sigma}_k{\bf F}^{-1}{\bf \Sigma}_k).
\end{equation}
By applying the properties of partial ordering, it follows %%%%%%%%%%%%%%%%%%
\begin{equation}\notag
\begin{aligned}
{\bf G}^T{\bf G} \preceq \sigma_{k+1}^4\|\widetilde{\bf \Omega}_{2}\widetilde{\bf \Omega}_{1}^{-1}\|_2^2 {\bf \Sigma}_k^{-4}  = {\bf \Psi}^4\|\widetilde{\bf \Omega}_{2}\widetilde{\bf \Omega}_{1}^{-1}\|_2^2, 
\end{aligned}
\end{equation}
where ${\bf \Psi} = \text{diag}(\psi_1, ..., \psi_k)$ is a $k \times k$ matrix with entries $\psi_i= \frac{\sigma_{k+1}}{\sigma_i}$. Moreover %%%%%%%%%%%%%%%%%%
\begin{equation}\notag
 {\bf \Sigma}_k({\bf I} + {\bf G}^T{\bf G})^{-1}{\bf \Sigma}_k \succeq {\bf \Sigma}_k({\bf I} + {\bf \Psi}^4\|\widetilde{\bf \Omega}_{2}\widetilde{\bf \Omega}_{1}^{-1}\|_2^2)^{-1}{\bf \Sigma}_k,
\end{equation}
which results in %%%%%%%%%%%%%%%%%%%%
\begin{equation}\notag
\begin{aligned}
\lambda_i({\bf A}^T{\bf A}) & \ge \lambda_i({\bf A}^T\widetilde{\bf Q}\widetilde{\bf Q}^T{\bf A})  \ge \lambda_i({\bf \Sigma}_k({\bf I} + {\bf G}^T{\bf G})^{-1}{\bf \Sigma}_k) \\ & \ge \frac{\sigma_i^2}{1 + \psi_i^4\|\widetilde{\bf \Omega}_{2}\widetilde{\bf \Omega}_{1}^{-1}\|_2^2}.
\end{aligned}
\end{equation}
The desired result follows by taking the square root of the last identity.

%%%%%%%%%%%%%%%%%%%%        %%%%%%%%%%%%%%%%%%%%        %%%%%%%%%%%%%%%%%%

\subsubsection{Proof of equation \eqref{eqL22B}}
Let 
\begin{equation}\notag
\bar{\bf K} \triangleq {\bf A}^T{\bf Q} = [\bar{\bf K}_1\quad \bar{\bf K}_2]=[{\bf P}_1\quad {\bf P}_2]
\begin{bmatrix}
{\bf L}_{11}^T & {\bf L}_{21}^T  \\
{\bf 0} & {\bf L}_{22}^T
\end{bmatrix},
\end{equation}
where $\bar{\bf K}_1$ and ${\bf P}_1$ contain the first $k$ columns, and $\bar{\bf K}_2$ and ${\bf P}_2$ contain the remaining $n-k$ columns of $\bar{\bf K}$ and ${\bf P}$, respectively. It follows that 

\begin{equation}\notag
\begin{aligned}
\bar{\bf K}_1 = {\bf P}_1{\bf L}_{11}^T, \quad  \text{and} \quad
\bar{\bf K}_2 = {\bf V}_k{\bf \Sigma}_k{\bf U}_k^T{\bf Q}_2 + {\bf V}_{\perp}{\bf \Sigma}_{\perp}{\bf U}_\perp^T{\bf Q}_2.
\end{aligned}
\end{equation}
We then obtain (see, e.g., \cite[p. 13]{StewartSun90}) the following relation for ${\bf L}_{22}$:

\begin{equation}\notag
\begin{aligned}
{\bf P}_{\bar{\bf K}_1^\perp} \bar{\bf K}_2 = {\bf P}_2{\bf L}_{22}^T = ({\bf I} - \bar{\bf K}_1\bar{\bf K}_1^{-1})\bar{\bf K}_2  = ({\bf I} - {\bf P}_1{\bf P}_1^T)\bar{\bf K}_2.
\end{aligned}
\end{equation}
Hence %%%%%%%%%%%%%%%%%%%%%%%   
\begin{equation}\label{eqL22BoTwoTerms}
\begin{aligned}
\|{\bf L}_{22} \|_2 & \le \|{\bf P}_2^T\|_2\|({\bf I} - {\bf P}_1{\bf P}_1^T)\bar{\bf K}_2\|_2
\le \|({\bf I} - {\bf P}_1{\bf P}_1^T)\bar{\bf K}_2\|_2\\
& \le \|({\bf I} - {\bf P}_1{\bf P}_1^T){\bf A}_k^T\|_2 +\|({\bf I} - {\bf P}_1{\bf P}_1^T){\bf V}_\perp{\bf \Sigma}_{\perp}\|_2,
\end{aligned}
\end{equation}
The first line is due to the orthonormality of ${\bf P}_2$, while the second line follows because of the orthonormality of ${\bf Q}_2$ as well as the triangle inequality, where ${\bf A}_k^T \triangleq {\bf V}_k{\bf \Sigma}_k{\bf U}_k^T$. For the second term on the right-hand side of \eqref{eqL22BoTwoTerms}, we get
\begin{equation} \label{eq_2ndTermL22}
\begin{aligned}
\|({\bf I} - {\bf P}_1{\bf P}_1^T){\bf V}_\perp{\bf \Sigma}_{\perp}\|_2 \le \|({\bf I} - {\bf P}_1{\bf P}_1^T)\|_2\|{\bf V}_\perp\|_2\|{\bf \Sigma}_{\perp}\|_2 \le \sigma_{k+1}.
\end{aligned}
\end{equation}
%%%%%%%%%%%%%%%%%%%%%%%%%%%%%

To bound the first term on the right-hand side of \eqref{eqL22BoTwoTerms}, we need some preparation. We have for the range of the matrix $\bf P$ in Step 5 of Algorithm 1:
\begin{equation}\notag
\mathcal{R}({\bf P}) = \mathcal{R}({\bf A}^T{\bf Q})= \mathcal{R}({\bf A}^T{\bf A}\bar{\bf Q})= \mathcal{R}({\bf A}^T{\bf AA}^T {\bf \Omega}).
\end{equation}
Thus
\begin{equation}\notag
\mathcal{R}({\bf P}_1) = \mathcal{R}({\bf A}^T{\bf Q}_1)= \mathcal{R}({\bf A}^T{\bf A}\bar{\bf Q}_1)= \mathcal{R}({\bf A}^T{\bf AA}^T {\bf \Omega}_1).
\end{equation}
Defining a $k\times k$ matrix $\bar{\bf X} \triangleq \widetilde{\bf \Omega}_{1}^{-1}{\bf \Sigma}_k^{-3}$ and performing a QR factorization on ${\bf A}^T{\bf AA}^T {\bf \Omega}_1\bar{\bf X}$ gives 
\begin{equation} \label{eqP1eqQring}
\begin{aligned}
{\bf A}^T{\bf AA}^T {\bf \Omega}_1\bar{\bf X} = {\bf V}
\begin{bmatrix}
{\bf I}_k  \\
{\bf \Sigma}_\perp^3\widetilde{\bf \Omega}_{2}\widetilde{\bf \Omega}_{1}^{-1} {\bf \Sigma}_k^{-3} 
\end{bmatrix}= \mathring{\bf Q}\mathring{\bf R}.
\end{aligned}
\end{equation}
Let $\bar{\bf G} = {\bf \Sigma}_\perp^3\widetilde{\bf \Omega}_{2}\widetilde{\bf \Omega}_{1}^{-1} {\bf \Sigma}_k^{-3}$, and 
\begin{equation}\notag
\begin{aligned}
\bar{\bf F}^{-1} \triangleq \mathring{\bf R}^{-1} \mathring{\bf R}^{-T} = 
(\mathring{\bf R}^T \mathring{\bf R})^{-1} = ({\bf I}+\bar{\bf G}^T\bar{\bf G}) ^{-1}.
\end{aligned}
\end{equation}
Thus
\begin{equation}\notag
\mathring{\bf Q}\mathring{\bf Q}^T = {\bf V}
\begin{bmatrix}
\bar{\bf F}^{-1} & \bar{\bf F}^{-1}\bar{\bf G}^T  \\
\bar{\bf G}\bar{\bf F}^{-1} & \bar{\bf G}\bar{\bf F}^{-1}\bar{\bf G}^T 
\end{bmatrix}{\bf V}^T.
\end{equation}
%%%%%%%%%%%%%%%%%%%%%%%%%%%%
Since ${\bf P}_1{\bf P}_1^T = \mathring{\bf Q}\mathring{\bf Q}^T$, we get
\begin{equation}\notag
{\bf I} - {\bf P}_1{\bf P}_1^T = {\bf I} - \mathring{\bf Q}\mathring{\bf Q}^T = {\bf V}
\begin{bmatrix}
{\bf I} - \bar{\bf F}^{-1} & \bar{\bf F}^{-1}{\bf G}^T  \\
-\bar{\bf G}\bar{\bf F}^{-1} & {\bf I} -\bar{\bf G}\bar{\bf F}^{-1}\bar{\bf G}^T 
\end{bmatrix}{\bf V}^T.
\end{equation}
Writing ${\bf A}_k^T = {\bf V}[{\bf \Sigma}_k \quad {\bf 0}; {\bf 0} \quad {\bf 0}]{\bf U}^T$, we obtain
\begin{equation}\notag
\begin{aligned}
({\bf I} - {\bf P}_1{\bf P}_1^T){\bf A}_k^T = {\bf V}
	\begin{bmatrix}
	({\bf I} - \bar{\bf F}^{-1}){\bf \Sigma}_k  \\
	- \bar{\bf G} \bar{\bf F}^{-1}{\bf \Sigma}_k
	\end{bmatrix} {\bf U}^T.
\end{aligned}
\end{equation}
It follows that
\begin{equation}
\begin{aligned}\label{eqSigmKGbTGbSigmK}
\|({\bf I} - {\bf P}_1{\bf P}_1^T){\bf A}_k^T \|_2^2  = \|{\bf \Sigma}_k ({\bf I} - \bar{\bf F}^{-1}){\bf \Sigma}_k\|_2 \le \|{\bf \Sigma}_k \bar{\bf G}^T\bar{\bf G}{\bf \Sigma}_k \|_2.
\end{aligned} 
\end{equation}
To get the equality, we have used the following relation that holds for any matrix $\bf A$ with ${\bf I} +{\bf A}$ being non-singular \cite{HendersonS81}:
\begin{equation}\label{eqInvEqs}
({\bf I} +{\bf A})^{-1} = {\bf I} - {\bf A}({\bf I} +{\bf A})^{-1} = {\bf I} - ({\bf I} +{\bf A})^{-1}{\bf A}.
\end{equation}
The inequality in \eqref{eqSigmKGbTGbSigmK} follows because
${\bf I} - \bar{\bf F}^{-1}=\bar{\bf G}^T\bar{\bf G}({\bf I}+\bar{\bf G}^T\bar{\bf G}) ^{-1}\preccurlyeq \bar{\bf G}^T\bar{\bf G}$. Therefore

\begin{equation}
\begin{aligned}\notag
\|({\bf I} - {\bf P}_1{\bf P}_1^T){\bf A}_k^T \|_2  \le \|\bar{\bf G}{\bf \Sigma}_k \|_2 \le \psi_k^2\|\widetilde{\bf \Omega}_{2}\widetilde{\bf \Omega}_{1}^{-1} \|_2\sigma_{k+1}.
\end{aligned} 
\end{equation}

Inserting this result and that of \eqref{eq_2ndTermL22} into \eqref{eqL22BoTwoTerms} gives the desired bound.   %\QEDB

%%%%%%%%%%%%%%%%%%%%%%%%%%%%%%%%++++++++++++++++++++++++++++++++++++++++++++
%%%%%%%%%%%%%%%%%%%%%%%%%%%%%%%%++++++++++++++++++++++++++++++++++++++++++++

\subsection{Proof of Theorem \ref{ThSins}}
\label{AppdxPrTh2}

{\textit{Upper bound for $\text{sin}\theta_Q$.} According to \cite[Theorem 2.6.1]{GolubVanLoan96} we have
\begin{equation}\notag
\begin{aligned}
\text{sin}\theta_Q & = \|{\bf U}_\perp^T{\bf Q}_1\|_2 = \|{\bf U}_k{\bf U}_k^T - {\bf Q}_1{\bf Q}_1^T\|_2 \\
& = \|{\bf U}_k{\bf U}_k^T - \widetilde{\bf Q}\widetilde{\bf Q}^T\|_2 = \|{\bf U}_\perp^T \widetilde{\bf Q}\|_2. 
%\label{eq_dist_def}
\end{aligned}
\end{equation}
The first relation in the second line follows due to ${\bf Q}_1{\bf Q}_1^T = \widetilde{\bf Q}\widetilde{\bf Q}^T$, as shown before. Writing \eqref{eq_QRtilde} as:
\begin{equation}\notag
\begin{aligned}
\begin{bmatrix}
{\bf I}_k  \\
{\bf \Sigma}_\perp^2\widetilde{\bf \Omega}_{2}\widetilde{\bf \Omega}_{1}^{-1} {\bf \Sigma}_k^{-2}
\end{bmatrix} = \begin{bmatrix}
{\bf U}_k^T  \\
{\bf U}_\perp^T 
\end{bmatrix}\widetilde{\bf Q}\widetilde{\bf R},
\end{aligned}
\end{equation}
it follows $\widetilde{\bf R}^{-1} = {\bf U}_k^T \widetilde{\bf Q}$, and hence  %%%%%%%%%%%%
\begin{equation}\notag
\begin{aligned}
{\bf U}_\perp^T \widetilde{\bf Q} = {\bf \Sigma}_\perp^2\widetilde{\bf \Omega}_{2}\widetilde{\bf \Omega}_{1}^{-1} {\bf \Sigma}_k^{-2}{\bf U}_k \widetilde{\bf Q}.
\end{aligned}
\end{equation}
Therefore %%%%%%%%
\begin{equation}\notag
\begin{aligned}
\text{sin}\theta_Q = 
\|{\bf U}_\perp^T \widetilde{\bf Q}\|_2  \le \|{\bf \Sigma}_\perp^2\widetilde{\bf \Omega}_{2}\widetilde{\bf \Omega}_{1}^{-1} {\bf \Sigma}_k^{-2}\|_2\|{\bf U}_k \widetilde{\bf Q}\|_2 \le \psi^2_k\|{\bf \Omega}_{2}{\bf \Omega}_{1}^{-1}\|_2.
\end{aligned}
\end{equation}

%%%%%%%%%%%%%%%%%%%%%%%%%%%%%%%%%%%%%%%%%%%%%%%%%%%% 

\textit{Upper bound for $\text{sin}\theta_P$}. $\text{sin}\theta_P$ is defined as follows: 
\begin{equation}\notag
\begin{aligned}
\text{sin}\theta_P & = \|{\bf V}_\perp^T {\bf P}_1\|_2  =  \|{\bf V}_k{\bf V}_k^T - {\bf P}_1{\bf P}_1^T\|_2 \\ & = \|{\bf V}_k{\bf V}_k^T - \mathring{\bf Q}\mathring{\bf Q}^T\|_2 =  \|{\bf V}_\perp^T \mathring{\bf Q}\|_2.
\end{aligned}
\end{equation}
The first relation in the second line follows because ${\bf P}_1{\bf P}_1^T = \mathring{\bf Q}\mathring{\bf Q}^T$. From \eqref{eqP1eqQring}, we get $\mathring{\bf R}^{-1} = {\bf V}_k^T \mathring{\bf Q}$. Henec
\begin{equation}\notag
\begin{aligned}
\text{sin}\theta_P = \|{\bf V}_\perp^T \mathring{\bf Q}\|_2 =\| {\bf \Sigma}_\perp^3\widetilde{\bf \Omega}_{2}\widetilde{\bf \Omega}_{1}^{-1} {\bf \Sigma}_k^{-3}{\bf V}_k \mathring{\bf Q}\|_2 \le \psi^3_k\|{\bf \Omega}_{2}{\bf \Omega}_{1}^{-1}\|_2.
\end{aligned}
\end{equation} 

%%%%%%%%%%%%%%%%%%%%%%%%%%%%%%%%%%%%%%%%%%%%%

\textit{Upper bound for $\text{sin}\phi_Q$}. 

\begin{equation}\notag
\begin{aligned}
\text{sin}\phi_Q & = \|{\bf U}_k^T{\bf Q}_2\|_2 =  \|{\bf U}_k^T{\bf Q}_2{\bf Q}_2\|_2 = \|{\bf U}_k^T({\bf I}-{\bf Q}_1{\bf Q}_1^T)\|_2 \\
& = \|{\bf U}_k^T({\bf I}-\widetilde{\bf Q}\widetilde{\bf Q}^T)\|_2 
  = \|[\underbrace{{\bf I}-{\bf F}^{-1} \quad {\bf F}^{-1}{\bf G}^T]{\bf U}^T}_{ \triangleq \widetilde{\bf M}} \|_2.
\end{aligned}
\end{equation}
The first equality follows due to the orthonormality of ${\bf Q}_2$, and the third equality from the proof of Theorem \ref{ThSinRQLP}, equation \ref{eqQ1eqQtilde}.  It follows that    
\begin{equation}
\begin{aligned}\notag
\widetilde{\bf M} \widetilde{\bf M}^T = {\bf I} - {\bf F}^{-1},
\end{aligned} 
\end{equation}
where we have used equation \eqref{eqInvEqs}. Since the matrix ${\bf I}-{\bf F}^{-1} = {\bf G}^T{\bf G}({\bf I}+{\bf G}^T{\bf G})^{-1}$ is positive semidefinite, its eigenvalues satisfy \cite[p. 148]{StewartSun90}:

\begin{equation}
\begin{aligned}
\lambda_i({\bf I}-{\bf F}^{-1}) = \frac{{\sigma}_i^2({\bf G})}{1+{\sigma}_i^2({\bf G})}, \quad i=1,...k.
\end{aligned}
\notag
\end{equation}
The largest singular value of ${\bf G}$ satisfies:
\begin{equation}\notag
 \sigma_1 ({\bf G}) \le \psi_k^2\|\widetilde{\bf \Omega}_{2}\widetilde{\bf \Omega}_{1}^{-1} \|_2.
\end{equation}
It therefore follows that
\begin{equation}
\begin{aligned}
\text{sin}^2\phi_Q =  \lambda_1({\bf I}-{\bf F}^{-1}) \le \frac{\psi_k^4\|\widetilde{\bf \Omega}_{21}\widetilde{\bf \Omega}_{11}^{-1} \|_2^2}{1+\psi_k^4\|\widetilde{\bf \Omega}_{2}\widetilde{\bf \Omega}_{1}^{-1} \|_2^2}.
\end{aligned} \notag
\end{equation}
The desired result is obtained  by taking the square root.
%%%%%%%%%%%%%%%%%%%%%%%%%%%%%%%%%%%%%%%%%%%%%%

\textit{Upper bound for $\text{sin}\phi_P$.} 
\begin{equation}\notag
\begin{aligned}
\text{sin}\phi_P & = \|{\bf V}_k^T{\bf P}_2\|_2 =  \|{\bf V}_k^T{\bf P}_2{\bf P}_2\|_2 = \|{\bf V}_k^T({\bf I}-{\bf P}_1{\bf P}_1^T)\|_2 \\
& = \|{\bf V}_k^T({\bf I}-\mathring{\bf Q}\mathring{\bf Q}^T)\|_2 
  = \|[{\bf I}-\bar{\bf F}^{-1} \quad \bar{\bf F}^{-1}\bar{\bf G}^T]{\bf V}^T \|_2.
\end{aligned}
\end{equation}
Following the same procedure as the proof of the bound for $\text{sin}\phi_Q$, the result follows.

\subsection{Proof of Theorem \ref{ThLowRankErr}}
\label{AppdxPrTh3}
Using ${\bf A}={\bf A}_k+{\bf A}_\perp$, we will obtain
\begin{equation}\label{eqProoLREE}
\begin{aligned}
\|({\bf I} - {\bf Q}_k{\bf Q}_k^T){\bf A}\|_\rho & \le \|({\bf I} -{\bf Q}_k{\bf Q}_k^T){\bf A}_k\|_\rho + \|({\bf I} - {\bf Q}_k{\bf Q}_k^T){\bf A}_\perp\|_\rho \\
& \le \|({\bf I} - {\bf Q}_1{\bf Q}_1^T){\bf A}_k\|_\rho + \|{\bf \Sigma}_\perp\|_\rho.
\end{aligned}
\end{equation}
The first relation is due to the triangle inequality, while the second term on the right-hand side of the second relation is due to the strong submultiplicativity property \cite[equation 9.3.13]{Bernstein09}, which holds for any matrices $\bf A$ and $\bf B$ with conformal dimensions:   
\begin{equation}\notag
\begin{aligned}
  \|{\bf AB}\|_F \le \|{\bf A}\|_F\|{\bf B}\|_2,  \quad \text{and} \quad
 \|{\bf AB}\|_F \le \|{\bf A}\|_2\|{\bf B}\|_F.    
\end{aligned}
\end{equation}
Then
\begin{equation}\notag
\begin{aligned}
 \|({\bf I} - {\bf Q}_1{\bf Q}_1^T){\bf A}_k\|_\rho & = \|({\bf I} - {\bf Q}_1{\bf Q}_1^T){\bf U}_k{\bf U}_k^T{\bf A}\|_\rho  \le \text{sin}\theta_Q \|{\bf A}\|_\rho.
\end{aligned}
\end{equation}
Inserting the above result into \eqref{eqProoLREE} gives the desired bound. The bound for $\|{\bf A}({\bf I} - {\bf P}_k{\bf P}_k^T)\|_\rho$ follows similarly.

\bibliography{mybibfile}

\end{document}